\documentclass[12pt]{article}
\usepackage[utf8]{inputenc}
\usepackage[french,english]{babel}
\usepackage[T1]{fontenc}
\usepackage{amsmath}
\usepackage{amsfonts}
\usepackage{amssymb}
\usepackage{graphicx}
\usepackage[dvipsnames]{xcolor}
\usepackage{tcolorbox}
\usepackage{hyperref}
\usepackage[nottoc]{tocbibind}
\usepackage[right=2cm, left=2cm, top=2cm, bottom=2cm]{geometry}
\usepackage{cite}

\newcommand{\overbar}[1]{\mkern 1.5mu\overline{\mkern-1.5mu#1\mkern-1.5mu}\mkern 1.5mu}

\newtheorem{Theorem}{Theorem}[section]
\newtheorem{Lemma}{Lemma}[section]
\newtheorem{Proposition}{Proposition}[section]
\newtheorem{Definition}{Definition}[section]
\newtheorem{Remark}{Remark}[section]
\newtheorem{Corollary}{Corollary}[section]

\newgeometry{includefoot,left=2cm,right=2cm,bottom=1cm,top=2cm}
\title{ \LARGE \bf%
{$L^p$-asymptotic stability of 1D damped wave equations with localized and linear damping}
\thanks{This research was partially supported by the iCODE Institute, research project of the IDEX Paris-Saclay, and by the Hadamard Mathematics LabEx (LMH) through the grant number ANR-11-LABX-0056-LMH in the ``Programme des Investissements d'Avenir''.}
}
\author{Meryem Kafnemer$^{1}$, Benmiloud Mebkhout$^{2}$, Fr\'{e}d\'{e}ric Jean$^{1}$, and Yacine Chitour$^{3}$}
\footnotetext[1]{UMA, ENSTA Paris, Institut Polytechnique de Paris, F-91120 Palaiseau, France}
\footnotetext[2]{D\'epartement de Math\'ematiques, Universit\'{e} Abou Bekr Belkaid, Tlemcen, Algeria.}
\footnotetext[3]{L2S, Universit\'e Paris Saclay, France.}
\begin{document}
\maketitle
\begin{abstract}
In this paper, we study the $L^p$-asymptotic stability of the one dimensional linear damped
wave equation with Dirichlet boundary conditions in $[0,1]$, with $p\in (1,\infty)$. The damping
term is assumed to be linear and localized  to an arbitrary open sub-interval of $[0,1]$. We prove that the 
semi-group $(S_p(t))_{t\geq 0}$ associated with the previous equation is well-posed and exponentially stable.
The proof relies on the multiplier method and depends on whether $p\geq 2$ or $1<p<2$.
\end{abstract}
\section{Introduction}
This paper is concerned with the asymptotic stability of the one dimensional wave equation with a localized damping term and Dirichlet boundary conditions. The problem  is written as follows
\begin{equation}\label{prob}
\left\lbrace
\begin{array}{cccc}
z_{tt} -z_{xx} + a(x)z_t=0 & \hbox{for } (t,x) \in  \mathbb{R_+} \times (0,1), \\
z(t,0)= z(t,1)=0   &  t\geq 0 ,\\
z(0,\cdot)= z_0\ ,\ z_t(0,\cdot)=z_1,
\end{array}
\right.
\end{equation}
where $z$ is the solution of the problem, $(z_0,z_1)$ are the initial conditions and they all belong to an 
$L^p$-based functional space that will be defined later. The function $a$ is a continuous non-negative 
function on $[0,1]$, bounded from below by a positive constant on some non-empty open interval $\omega$ 
of $(0,1)$, which represents the region of the domain where the damping term is active.

Problem \eqref{prob} has been widely studied in the case $p=2$ whether with a linear or a non-linear 
damping. Stability results are proved under a geometric condition imposed on the damping domain 
$\omega$: it is properly introduced in the early work \cite{Zuazua1990} where the semi-linear problem is 
considered even in higher dimension and the geometric condition is extended and characterized in 
\cite{Liu1997}. Moreover, for linear problems there exist necessary and sufficient geometrical conditions for 
stabilization based on geometric optics methods (cf.\ the seminal work \cite{Bardos1992}). Strong 
stabilization, i.e., energy decay to zero for each trajectory, has been established in \cite{Dafermos} and 
\cite{Haraux} with a LaSalle's invariance argument. For the linear localized damping case in higher 
dimensions, exponential stability has been established several times using different tools, in particular using 
the multiplier method which is the relevant method to our paper context. We refer the reader to 
\cite{Komornik1994} for a complete presentation of the method as well as the tools associated to it. As for the 
stability results obtained by this method in this case, we refer for instance to \cite{Alabau} and 
\cite{Martinez1999} for detailed proofs and extended references. The non-linear problem on the other hand 
has been studied (for instance) in  \cite{Martinez2000} with no localization and in \cite{Kafnemer2020} for a 
localized damping. We refer the reader to the excellent survey \cite{Alabau} for more references in the 
Hilbertian framework i.e.\ when $p=2$.

As for more general functional frameworks, in particular $L^p$-based spaces with $p\neq 2$, few results exist 
and one reason is probably due to the fact that, in such non-Hilbertian framework, the semi-group associated 
with the d'Alembertian (i.e., the linear operator defining the wave equation) is not defined in general as soon 
as the space dimension is larger than or equal to two, see e.g., \cite{Peral1980}). This is why most of the 
existing results focus on several stabilization issues only in one spatial dimension. Well-posedness results as 
well as important $L^p$ estimates have been shown in \cite{Haraux}, in particular the introduction of a $p$-th 
energy of a solution as a generalization of the standard $E_2(t)=\int_0^1\frac{z_x^2+z_t^2}2$. Some of these 
results have been used in \cite{Amadori2020,Chitour2019} recently. The latter reference relies on Lyapunov 
techniques for linear time varying systems to prove $L^p$ exponential stability in the nonlinear problem under 
the hypothesis that initial data live in $L^\infty$ functional spaces and with $ p\geq 2$ only; other stability 
results have been shown in the same reference in particular $L^\infty$ stability but always with more 
conditions on initial data which creates a difference between the norms of trajectories and the norms of initial 
data used in their decay estimates.

In this paper we extend the results existing in the case $p=2$ to the case $p\in (1,\infty)$ by adapting the multiplier method to that issue. We start first by stating the problem and defining the appropriate $L^p$ functional framework as well as the notion of solutions. We prove the well-posedness of the corresponding $C^0$ semi-group of solutions using an argument inspired by \cite{Haraux2009} and \cite{Chitour2019}. As for stability issue, we prove that these semi-groups are indeed exponentially stable. Even though the argument depends on whether $p\geq 2$ or $p\in (1,2)$, it is another instance of the multiplier method, where the multipliers are expressed in terms of the Riemann invariants coordinates $\rho=z_x+z_t$ and $\xi=z_x-z_t$. In particular, one of the multipliers in the case $p=2$ is equal to $\phi(x)z$ with $\phi$ a non negative function which is used to localize estimates inside $\omega$. If $p\geq 2$, this multiplier is replaced by the pair of functions $\phi(x)z\vert \rho\vert^{p-2}$ and $\phi(x)z\vert \xi\vert^{p-2}$. Clearly, such multipliers cannot be used directly if $p\in (1,2)$ and must be modified, which yields to a more delicate treatment. In both cases, energy integral estimates are established following the standard strategy of the multiplier method and exponential stability is proved. For the two extremes cases $p=1$ and $p=\infty$, we are able to prove that the corresponding semi-groups are exponentially stable only for particular cases of global constant damping. However, we conjecture that such a fact should be true in case of any localized damping.

 The paper is divided into four sections, the first one being the introduction and the second one devoted to provide the main notations used throughout the paper. Section $3$ deals with the well-posedness issue
 and Section $4$ contains the main result of the paper, i.e.\ exponential stability of the $C^0$ semi-group of solutions for $p\in (1,\infty)$ as well as the partial result for $p=1$ and $p=\infty$. We gather in an appendix several technical results.

\textbf{Acknowledgment: }
We would like to thank Dario Prandi, Cyprien Tamekue and Nicolas Lerner for helpful discussions.

\section{Statement and main notations of the problem}

Consider Problem~\eqref{prob} where we assume the following hypothesis satisfied:\\ \\
$\mathbf{ (H_1)}$ $ a : [0,1] \rightarrow \mathbb{R}\ $ is a non-negative continuous function such that
\begin{equation} \label{a0}
\exists \ a_0 >0, \ a \geq a_0\ \ \hbox{on}\ \ \omega =[c,d] \subset [0,1],
\end{equation}
where $\omega$ is a non empty interval such that $c=0$ or $d=1$, i.e., $\bar{\omega}$ contains a
neighborhood of $0$ or $1$. There is no loss of generality in assuming $d=1$, taking $0$ as an observation point.
\begin{Remark}
The results of this paper still hold if the assumption that $c=0$ or $d=1$ is removed
by using a piecewise multiplier method, i.e., we can use both $0$ and $1$ as observation points
(instead of simply $0$ here) to obtain the required energy estimate.
\end{Remark}
For $p \in [1, \infty)$, consider the function spaces
\begin{align}
X_p&:= W^{1,p}_0(0,1) \times L^p(0,1), \\
Y_p&:= \left(W^{2,p}(0,1) \cap W^{1,p}_0(0,1)\right) \times W^{1,p}_0(0,1),
\end{align}
where $X_p$ is equipped with the norm
\begin{align}
\Vert{(u,v)}\Vert_{X_p}:= \left(\frac{1}{p}\int_0^1 \left(|u' + v|^p + |u' - v|^p\right) dx \right)^\frac{1}{p},
\end{align}
and the space $Y_p$ is equipped with the norm
\begin{align}
\Vert(u,v)\Vert_{Y_p}:=
\left(\frac{1}{p}\int_0^1 \left(|u'' + v'|^p + |u'' - v'|^p\right) dx \right)^\frac{1}{p}.
\end{align}
Initial conditions $(z_0,z_1)$ for weak (resp.\ strong) solutions of \eqref{prob} are taken in $ X_p$ (resp.
in $ Y_p$), where the two concepts of solutions are precisely defined later in
Definition~\ref{def:weak-strong}.

For all $(t,x) \in \mathbb{R_+}\times (0,1)$, define the Riemann invariants
\begin{align}
\rho(t,x)= z_x(t,x)+z_t(t,x), \\
\xi(t,x)= z_x(t,x)-z_t(t,x).
\end{align}
Along strong solutions of \eqref{prob}, we deduce that
\begin{equation}\label{probw}
\left\lbrace
\begin{array}{ll}
\rho_t - \rho_x  =-\frac{1}{2}a(x)(\rho-\xi) &\text{in }   \mathbb{R_+} \times (0,1), \\
\xi_t + \xi_x  =\frac{1}{2}a(x)(\rho-\xi) &\text{in }   \mathbb{R_+} \times (0,1), \\
\rho(t,0)-\xi(t,0)= \rho(t,1)-\xi(t,1)=0   &  \forall t \in \mathbb{R_+} ,\\
\rho_0:=\rho(0,.)= z_0'+ z_1\ ,\ \xi_0 :=\xi(0,.)=z_0'-z_1,
\end{array}
\right.
\end{equation}
with $\left(\rho_0, \xi_0 \right) \in W^{1,p}(0,1)\times W^{1,p}(0,1)$.

We define the $p$th-energy of a (weak) solution of \eqref{prob} as the function $E_p$ defined on
$\mathbb{R_+}$ by
\begin{align}
E_p(t)=\frac{1}{p}\int_0^1\left( |z_x + z_t|^p + |z_x - z_t|^p \right) dx
\end{align}
and $E_p$ can be expressed in terms of $\xi$ and $\rho$ as
\begin{align}
E_p(t)= \frac{1}{p} \int_0^1 (|{\rho}|^p+|{\xi}|^p) dx.
\end{align}
For $r\geq 0$, we introduce the following notation
\begin{align}
\lfloor x\rceil^r:=\textrm{sgn}(x)|x|^r ,\ \   \forall x \in \mathbb{R},
\end{align}
where $\textrm{sgn}(x)=\frac{x}{\vert x\vert}$ for nonzero $x\in\mathbb{R}$ and {$\textrm{sgn}(0)=[-1,1]$}.
We have the following obvious formulas which will be used all over the paper:
\begin{align}
\frac{d}{dx}(\lfloor x\rceil^r)=r|x|^{r-1}, \ \ \forall r\geq 1,\ x\in\mathbb{R},\\
\frac{d}{dx}(|x|^r)=r\lfloor x \rceil^{r-1}, \ \ \forall r> 1,\ x\in\mathbb{R}.
\end{align}
Before we state our results, we provide the following proposition
(essentially taken from \cite{Haraux2009}).
\begin{Proposition}\label{prop1}
Let $p \in [1, \infty)$ and suppose that a strong solution $z$ of \eqref{probw} is defined on
a non trivial interval $I \subset \mathbb{R_+}$ containing $0$, for some initial conditions $(z_0,z_1)\in
Y_p$. For $t\in I$, define
\begin{align}\label{phidef}
\Phi (t) := \int_0^1 [{\cal{F}}(\rho)+{\cal{F}}(\xi)] dx,
\end{align}
where ${\cal{F}} $ is a $C^1$ convex function. Then $\Phi$ is well defined
for $t\in I$ and satisfies
\begin{align}\label{phidec}
\frac{d}{dt} \Phi(t) =-\frac{1}{2}\int_0^1 a(x)(\rho - \xi)({\cal{F}}'(\rho)- {\cal{F}}'(\xi)) dx \leq 0.
\end{align}
\end{Proposition}
\textbf{\emph{Proof. }}
By the regularity assumptions, $\rho(t,.)$ and $\xi(t,.)$ are absolutely continuous functions. Formal
differentiation, easy to justify a posteriori by the regularity of the data, yields
\begin{align}\label{p11}
\frac{d}{dt} \int_0^1 [{\cal{F}}(\rho)+{\cal{F}}(\xi)] dx =
\int_0^1 [ \rho_t {\cal{F}}'(\rho) + \xi_t {\cal{F}}'(\xi)] dx.
\end{align}
Using \eqref{probw}, one obtains that
\begin{align}
\frac{d}{dt} \int_0^1 ({\cal{F}}(\rho)+{\cal{F}}(\xi)) dx
&= \int_0^1(\rho_x -\frac{1}{2} a(x)(\rho -\xi)) {\cal{F}}'(\rho) +
(-\xi_x +\frac{1}{2} a(x)(\rho -\xi)) {\cal{F}}'(\xi)) dx,\notag\\
&= \int_0^1 [{\cal{F}}(\rho)-{\cal{F}}(\xi)]_x\, dx -
\frac{1}{2}\int_0^1 a(x)(\rho - \xi)({\cal{F}}'(\rho)- {\cal{F}}'(\xi)) dx, \notag\\
&=- \frac{1}{2}\int_0^1 a(x)(\rho - \xi )({\cal{F}}'(\rho)- {\cal{F}}'(\xi)) dx. \label{p12}
\end{align}
Since ${\cal{F}}$ is convex, ${\cal{F}}'$ is non-decreasing, implying that
$(\rho - \xi)({\cal{F}}'(\rho)- {\cal{F}}'(\xi)) \geq 0$ which gives the conclusion when combining it with \eqref{p12}.
\begin{small}
\begin{flushright}
$\blacksquare$
\end{flushright}
\end{small}
\begin{Corollary}\label{cor:decrease}
For $(z_0,z_1) \in Y_p$, suppose that the solution $z$ of \eqref{prob}
exists on $\mathbb{R_+}$. Then the energy $t\longmapsto E_p(t)$ is non-increasing and, for $t\geq 0$,
\begin{align}\label{E_p'}
E_p'(t) = -\frac{1}{2}\int_0^1 a(x)(\rho-\xi)\left(\lfloor \rho \rceil^{p-1}-
\lfloor \xi \rceil^{p-1} \right) dx.
\end{align}
\end{Corollary}
\textbf{\emph{Proof. }}
For $(z_0,z_1) \in Y_p$ and $p>1$, we apply Proposition \ref{prop1} with
$F(s)=\frac{|s|^p}{p}$, which proves \eqref{E_p'}.
\begin{small}
\begin{flushright}
$\blacksquare$
\end{flushright}
\end{small}

\section{Well-posedness}

We start by recalling the classical representation formula for regular solutions of \eqref{prob}
given by the d'Alembert formula, cf. \cite[Equation 8, page 36]{Strauss2007}.

\begin{Proposition}\label{propdal}
Consider the following problem with an arbitrary source term
$g \in C^2(\mathbb{R_+} \times \mathbb{R},\mathbb{R})$
and initial data $z_0 \in C^2(\mathbb{R})$ and
$z_1 \in C^1(\mathbb{R})$,
\begin{equation}\label{probgR}
\left\lbrace
\begin{array}{cccc}
z_{tt}(t,x) -z_{xx}(t,x) + g(t,x)=0 & \text{for }  (t,x)\in \mathbb{R_+} \times \mathbb{R}, \\
z(0,.)= z_0\ ,\ z_t(0,.)=z_1.
\end{array}
\right.
\end{equation}
Then the unique solution $z$ of this problem is in
$C^2(\mathbb{R_+}\times \mathbb{R},\mathbb{R})$ and is given for all
$(t,x) \in \mathbb{R_+} \times \mathbb{R}$ by d'Alembert formula
\begin{align}\label{solR}
{z}(t,x)=\frac{1}{2}\left[{z_0}(x+t) +{z_0}(x-t)\right]+\frac{1}{2} \int_{x-t}^{x+t}{z_1}(s)ds+\frac{1}{2}\int_0^t\int_{x-(t-s)}^{x+(t-s)} {g}(s,\tau)\, d\tau\, ds.
\end{align}
\end{Proposition}
In order to apply the above proposition to \eqref{prob}, we extend by a standard procedure (cf. \cite[Exercise 4, section 4.3]{Folland1992}) the following partial differential
equation defined on $\mathbb{R_+} \times (0,1)$
\begin{equation}\label{probg}
\left\lbrace
\begin{array}{ll}
z_{tt} -z_{xx} + g(t,x)=0 & \text{for } (t,x) \in  \mathbb{R_+} \times (0,1), \\
z(t,0)= z(t,1)=0   &  \forall t \in \mathbb{R_+} ,\\
z(0,.)= z_0\ ,\ z_t(0,.)=z_1,
\end{array}
\right.
\end{equation}
to an equivalent partial differential system defined on
$\mathbb{R_+} \times \mathbb{R}$. We first extend the data of the problem by considering
$\tilde{z}_0$, $\tilde{z}_1$ and $\tilde{g}$ the $2$-periodic extensions to $\mathbb{R}$ of the odd
extensions  of $z_0$, $z_1$ and $g$ to $[-1,1]$.

Using \eqref{solR}, we obtain then the expression of the solution ${z}$ for the problem on $\mathbb{R_+} \times (0,1)$, which is the following
\begin{align}\label{z}
{z}(t,x)=\frac{1}{2}\left[\tilde{z}_0(x+t) +\tilde{z}_0(x-t)\right]+\frac{1}{2} \int_{x-t}^{x+t}\tilde{z}_1(s)ds+\frac{1}{2}\int_0^t\int_{x-(t-s)}^{x+(t-s)} \tilde{g}(s,\tau)\, d\tau\, ds,
\end{align}
which clearly provides, for every $t\geq 0$, a $2$-periodic odd function $z(t,\cdot)$ on $\mathbb{R}$.
We also have the expression of the derivatives
\begin{align}\label{zx}
z_x(t,x)=\frac{1}{2}\left[\tilde{z}_0'(x+t) +\tilde{z}_0'(x-t)\right]+\frac{1}{2} \left[\tilde{z}_1(x+t)-\tilde{z}_1(x-t)\right] \notag \\
+\frac{1}{2}\int_0^t \left[\tilde{g}(s,x+(t-s)) - \tilde{g}(s,x-(t-s))\right]\, ds,
\end{align}
and
\begin{align}\label{zt}
z_t(t,x)=\frac{1}{2}\left[\tilde{z}_0'(x+t) -\tilde{z}_0'(x-t)\right]+\frac{1}{2} \left[\tilde{z}_1(x+t)+\tilde{z}_1(x-t)\right] \notag \\
+\frac{1}{2}  \int_0^t \left[\tilde{g}(s,x+(t-s))+\tilde{g}(s,x-(t-s))\right] \,ds.
\end{align}
Before we proceed to the well-posedness of \eqref{prob} in $X_p$ (resp.\ $Y_p$), we need to define the
notion of its weak and strong solutions.
\begin{Definition}\label{def:weak-strong}
For $(z_0,z_1)$ in $X_p$ (resp.\ $Y_p$), we say that  \eqref{prob} has a weak (resp.\ strong)
solution  $z\in L^\infty(\mathbb{R_+},  W^{1,p}_0(0,1))\cap W^{1,\infty}(\mathbb{R_+},L^p(0,1))$
(resp.\ $z\in L^\infty(\mathbb{R_+}, W^{2,p}(0,1)\cap W^{1,p}_0(0,1))\cap W^{1,\infty}(\mathbb{R_+},
W^{1,p}_0(0,1))$) given by the expression \eqref{z} in $X_p$ (resp.\ $Y_p$), if the source term $g$ from
\eqref{probgR} is given by $g(t,x)= -\tilde{a}(x) z_t(t,x)$, where $\tilde{a}$ is the $2$-periodic extension
to $\mathbb{R}$ of the even extension of $a$ to $[-1,1]$.
\end{Definition}

\begin{Theorem}[Well-posedness]\label{thrmwp}
Let $p \in [1, \infty)$. For any initial data $(z_0,z_1)\in X_p$ (resp.\ $Y_p$), there exists a unique weak
(resp.\ strong) solution $z$ such that
\begin{align}
z \in L^\infty(\mathbb{R_+}, W^{1,p}_0(0,1))\cap W^{1,\infty}(\mathbb{R_+},L^p(0,1)),
\end{align}
\begin{align}\label{zztDp}
\hbox{(resp. } \qquad z\in L^\infty(\mathbb{R_+}, W^{2,p}(0,1)\cap W^{1,p}_0(0,1))\cap W^{1,\infty}(\mathbb{R_+}, W^{1,p}_0(0,1)).)
\end{align}
Moreover, in both cases, the energy function $t\mapsto E_p(t)$ associated with a solution is
non-increasing.
\end{Theorem}
\textbf{\emph{Proof.}} The arguments for both items is adapted from that of \cite[Theorem 1]
{Chitour2019}. We prove the existence of an appropriate solution $y$ of \eqref{zt} by a standard
fixed point argument. We proceed on some interval $[0,T]$ for $T>0$ small enough independent on
the initial condition. We can then reproduce the reasoning on $[T,2T]$ starting from the solution at $t=T$
and so on to establish well-posedness for all $t\geq 0$.

Since $\tilde{g}$ is $2$-periodic function in space, it is natural to work in a space of
functions that have the same features. Hence we denote by $\mathcal{B}_T$ the space of functions
that are defined on $\left[0,T\right] \times \mathbb{R}$, odd on $[-1,1]$ and
$2$-periodic in space and $p$-integrable. The space $\mathcal{B}_T$ is equipped with the norm
 \begin{align}
 \Vert y \Vert_{Y_T} = \sup_{t\in[0,T]} \Vert y(t,.) \Vert_{L^p(0,1)},
 \end{align}
which makes it a Banach Space. We define the mapping
\begin{align}
 F_T : \mathcal{B}_T &\longrightarrow \mathcal{B}_T \notag
 \\
 y &\longmapsto  F_T(y), \notag
 \end{align}
 such that, for all $(t,x) \in [0,T]\times \mathbb{R}$, we have
 \begin{align}
 &F_T(y)(t,x)=\frac{1}{2}\left[\tilde{z_0}'(x+t) -\tilde{z_0}'(x-t)\right]+\frac{1}{2} \left[\tilde{z_1}(x+t)+\tilde{z_1}(x-t)\right] \notag \\
 +\frac{1}{2} & \int_0^t \left[\tilde{a}(x+(t-s))y(s,x+(t-s))+\tilde{a}(x-(t-s))y(s,x-(t-s))\right] \,ds.
 \end{align}
 Since $a$ is bounded, it is clear that $F_T$ is a contraction on $\mathcal{B}_T$ for $T>0$ small
 enough, hence the existence of a fixed point to $F_T$, which is a (weak) solution of \eqref{prob}. It is
 also clear that $T$ does not depend on the initial condition $(z_0,z_1)\in X_p$. As explained
 previously, this enables one to prove well-posedness in $X_p$.

As for the part regarding $Y_p$, the argument is similar to the previous one, after replacing
$\mathcal{B}_T$ by the space $\mathcal{D}_T $ consisting of the functions defined on $\left[0,T\right] \times
\mathbb{R}$ which are odd on $[-1,1]$ and $2$-periodic in space with $p$-integrable derivative with respect to $x$,  equipped with the norm given by
\begin{align}
\Vert y \Vert_{\mathcal{D}_T} = \sup_{t\in[0,T]} \Vert y(t,.) \Vert_{W^{1,p}_0(0,1)}.
\end{align}
For $(z_0,z_1) \in X_p$, $p>1$, we get that $t\mapsto E_p(t)$ is non increasing by the fact that $Y_p$
is dense in $X_p$. For $p=1$, we use the facts that $X_p$ is dense in $X_1$ for $p>1$ and the map
$p\mapsto E_p(t)$, for a fixed trajectory and a fixed positive time $t$, is right-continuous.
\begin{small}
\begin{flushright}
$\blacksquare$
\end{flushright}
\end{small}
\begin{Remark}
Since \eqref{prob} is linear and $t\mapsto E_p(t)$ is non-increasing, the flow of its weak solutions
defines a $C^0$-semigroup $(S_p(t))_{t\geq 0}$ of contractions of $X_p$, for every $p\in[1,\infty)$.
\end{Remark}

\section{Exponential Stability}
In this section, we aim to establish exponential stability for the $C^0$-semigroup
$(S_p(t))_{t\geq 0}$ defining the weak solutions of \eqref{prob} for every $p\in(1,\infty)$. The argument
relies on the multiplier method and is slightly different whether $p\geq 2$ or not. Indeed, two multipliers
involve the exponent $p-2$, which becomes negative if $p\in (1,2)$. In the latter case, one must modify
all the multipliers to handle that situation.

Before starting describing such results, we have the following weaker general stability result for
$p\in(1,\infty)$.

\begin{Proposition}[Strong stability]\label{propstab} 
Fix $p \in [1, \infty)$ and suppose that Hypothesis $\mathbf{ (H_1)}$ is satisfied. Then, for every $(z_0,z_1)\in
X_p$, the solution $z(t,\cdot)$ of \eqref{prob} starting at $(z_0,z_1)$ tends to zero as $t$ tends to infinity.
\end{Proposition}
\textbf{\emph{Proof.}} We follow the proof provided in the case $p=2$ in \cite{Dafermos}: by a
standard density argument, it is enough to establish the result for strong solutions of \eqref{prob}. The
latter is obtained by a LaSalle type of argument using the energy function $E_p$ and the fact that
the set $\{z(t,\cdot),\ t\geq 0\}$ is relatively compact in $W^{1,p}_0(0,1)$, which is itself obtained by
noticing that $z_t$ is a weak solution of \eqref{prob} with bounded energy $E_p$. \hfill 
\begin{small}
\begin{flushright}
$\blacksquare$
\end{flushright}
\end{small}

We introduce next some functions and notations which are common to the handling of both cases.
Recalling that we have chosen
$x_0=0$ as an observation point, we consider for $0<\epsilon_0 < \epsilon_1 < \epsilon_2 $,
the sets $Q_i = (1-\varepsilon_i,1+\varepsilon_i )$,  $i=0,1,2$, as well as three smooth functions
$\psi,\phi$ and $\beta$ defined below according to the next figure.
\begin{center}
\includegraphics[width=15cm]{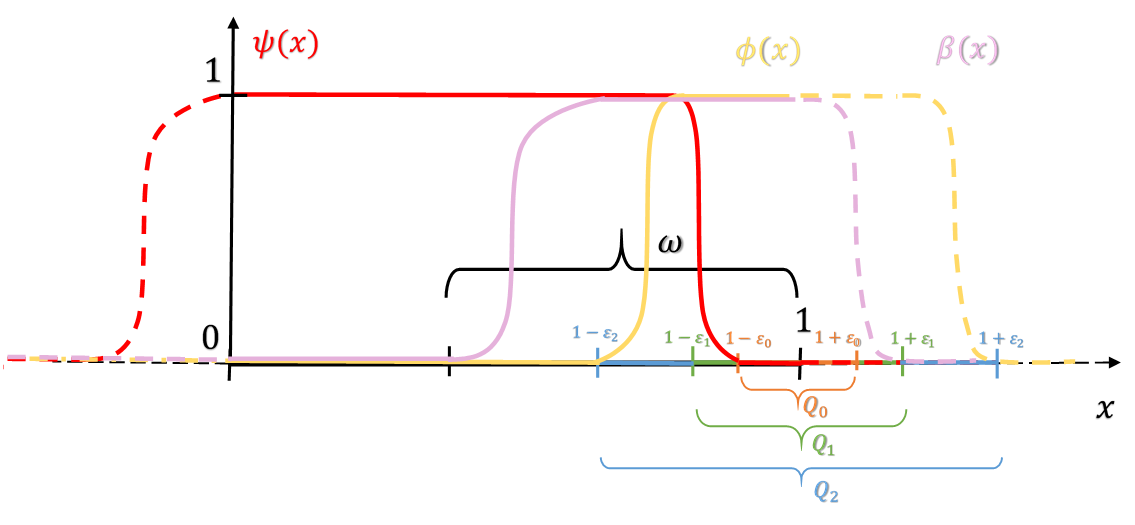}
\end{center}
More precisely, the functions $\psi,\phi$ and $\beta$ are smooth with compact support and
defined as follows:
\begin{equation}\label{psidef}
\begin{cases}
0\leq \psi \leq 1,
\\
\psi = 0\ \textrm{on}\ Q_0 ,
\\
\psi=1\ \textrm{on}\ (0,1)\setminus Q_1,
\end{cases}
\quad
\begin{cases}
0\leq \phi \leq 1,
\\
\phi = 1\ \hbox{on}\ Q_1 ,
\\
\phi=0\ \hbox{on}\ (0,1)\setminus Q_2,
\end{cases}
\quad \begin{cases}
0 \leq \beta \leq 1,
\\
\beta = 1\ \ \hbox{on}\ \ Q_2 \cap (0,1) ,
\\
\beta=0\ \ \hbox{on}\  \mathbb{R} \setminus \omega.
\end{cases}
\end{equation}

\begin{Remark}\label{rem:Cp}
In the sequel, we will denote by $C_p$ positive constants only depending on $p$ and by $C$ positive
constants depending on $a(\cdot)$ (typically through its upper bound $A$ on $[0,1]$ and its lower
bound $a_0$ on $\omega$), and on
 $\psi,\phi$ and $\beta$ (through bounds of their first derivatives over their supports).
\end{Remark}
 
Our main result is the following theorem:
 \begin{Theorem}\label{thrmstab} \textbf{(Exponential stability)}
Fix $ p \in (1, \infty)$ and suppose that Hypothesis $\mathbf{ (H_1)}$ is satisfied.  Then the $C^0$-semigroup
$(S_p(t))_{t\geq 0}$ is exponentially stable.
\end{Theorem}

\subsection{Case where $p\geq 2$}

As usual, it is enough to prove Theorem~\ref{thrmstab} for strong solutions and then extend the result for weak
solutions by a density argument. In turn, the theorem for strong solutions classically follows from the
next proposition, cf. \cite[Theorem 1.4.2]{Alabau} for instance.

\begin{Proposition}\label{prop}
Fix $ p \in [2, \infty)$ and suppose that Hypothesis $\mathbf{ (H_1)}$ is satisfied.  Then there exist positive constants $C$ and $C_p$
such that, for every $(z_0,z_1) \in Y_p$, it holds the following energy estimate:
\begin{equation} \label{expstab}
\forall \  0\leq S\leq T\ ,\  \int_S^T E_p(t)\, dt \leq C\,C_p E_p(S),
\end{equation}
where $E_p(\cdot)$ denotes the energy of the solution of \eqref{prob} starting at $(z_0,z_1)$.
 \end{Proposition}

The proof will be divided into four steps in subsections~\ref{se:411}--\ref{se:414}. We fix an arbitrary pair of times $0\leq S\leq T$ and a strong solution $z(\cdot,\cdot)$ of \eqref{prob} starting at $(z_0,z_1)\in Y_p$, and we consider three sets of multipliers:
\begin{description}
\item[$(m1)$] $x\mapsto x\psi(x) f(\rho(t,x))$ and $x\mapsto x\psi(x)f(\xi(t,x))$ for every $t\geq 0$;
\item[$(m2)$] $x\mapsto \phi(x) f'(\rho(t,x)) z(t,x)$ and $x\mapsto \phi(x) f'(\xi(t,x)) z(t,x)$
for every $t\geq 0$;
\item[$(m3)$] $x\mapsto v(t,x)$ for every $t\geq 0$, where $v$ is the solution of the following elliptic problem defined for every $t\geq 0$:
\begin{equation}\label{ellipt}
\left\lbrace
\begin{array}{ll}
v_{xx} = \beta f(z)  & x\in (0,1) , \\
v(0)=v(1)=0,&
\end{array}
\right.
\end{equation}
\end{description}
where the function $f$ is defined by
\begin{equation}\label{eq:f0}
f(s)=\lfloor s \rceil^{p-1},\qquad s\in\mathbb{R}.
\end{equation}
Introducing the function $F(s)= \int_0^s f(\tau) d\tau$, we have
\begin{equation}\label{eq:f01}
F(s) = \frac{|s|^p}{p},\quad F'=f,\qquad  f'(s)=(p-1)\vert s\vert^{p-2}.
\end{equation}

Note that we use the usual notation $q=\frac{p}{p-1}$ for the conjugate exponent of $p$.

\begin{Remark}\label{rem:p=2}
In the Hilbertian case $p=2$, the classical multipliers as given in \cite{Alabau} are
$x\psi(x)z_x(t,x)$, $x\phi(x)z(t,x)$ and $v$ associated with $p=2$ (i.e.\ $v_{xx} = \beta z$). Then, while clearly
our third multiplier $v$ is a straightforward extension of the Hilbertian case to any $p\in[1,\infty)$, the two
sets of multipliers given in Items $(m1)$ and $(m2)$ seem to be new, even if those of  Item $(m2)$
are identical when $p=2$.
\end{Remark}

\subsubsection{First set of multipliers}
\label{se:411}
The first step toward an energy estimate consists in obtaining an inequality that contains the expression
of the energy $E_p$ and, {for this purpose}, we use the first set of multipliers of Item $(m1)$.  We obtain the
following lemma.
\begin{Lemma}\label{lemma1m}
Under the hypotheses of Proposition~\ref{prop}, we have the following estimate
{\begin{align}\label{T}
\int_S^T E_p(t) dt \leq CC_p\,E_p(S) +
C\underbrace{\int_S^T \int_{Q_1\cap (0,1)} (F(\rho)+F(\xi)) \, dx\, dt}_{\mathbf{S_4}},
\end{align}}
where $C_{p}$ denotes constants that depend on $p$ only.
\end{Lemma}
\textbf{\emph{Proof.}}
Multiplying the first equation of \eqref{probw} by $x\, \psi\, f(\rho)$ and integrating over
$[S,T]\times [0,1]$, we obtain that
\begin{equation} \label{o1}
\int_S^T\int_0^1 x\, \psi\, f(\rho) \left(\rho_t - \rho_x
+\frac{1}{2}a(x)(\rho-\xi)\right)\, dx\, dt=0.
\end{equation}
Starting with $\int_S^T\int_0^1 x\, \psi\, f(\rho)\rho_t \, dx\, dt$, one has
\begin{align}\label{ipp1}
\int_S^T\int_0^1 x\, \psi\, f(\rho) \rho_t \, dx\, dt = \int_0^1 x\, \psi \int_S^T \, (F(\rho))_t dtdx = \int_0^1 x\, \psi  \, \left[F(\rho)\right]_S^T dx.
\end{align}
Regarding $-\int_S^T\int_0^1 x\, \psi\, f(\rho)\rho_x \, dx\, dt$, we use an integration by part with respect
to $x$ and we obtain
 \begin{align}\label{ipp2}
&-\int_S^T\int_0^1 x\, \psi\, f(\rho) \rho_x \, dx\, dt = -\int_S^T \int_0^1 x\, \psi  \, (F(\rho))_x \, dx\, dt
\notag \\
&= \int_S^T \int_0^1 (x\, \psi)_x  F(\rho) \, dx\, dt - \int_S^T \left[ x\, \psi  F(\rho) \right]_0^1 = \int_S^T
\int_0^1 (x\, \psi)_x  F(\rho) \, dx\, dt.
\end{align}
By combining \eqref{ipp1} and \eqref{ipp2}, we get
\begin{equation}\label{ip1}
 \int_S^T \int_0^1 (x\, \psi)_x  F(\rho) \, dx\, dt+ \int_0^1 x\, \psi \, \left[F(\rho)\right]_S^T dx
=-\frac12\int_S^T \int_0^1 x\psi(x) a(x) f(\rho)(\rho-\xi))\, dx\, dt.
\end{equation}
We proceed similarly by multiplying the second equation of \eqref{probw} by $x\, \psi\, f(\xi)$ and,
following the same steps that yielded \eqref{ip1}, we obtain that
\begin{align}\label{ip2}
\int_S^T \int_0^1 (x\, \psi)_x  F(\xi) \, dx\, dt- \int_0^1 x\, \psi \, \left[F(\xi)\right]_S^T dx =-\frac{1}{2}\int_S^T \int_0^1 x\psi(x) a(x) f(\xi)(\rho-\xi))\, dx\, dt.
\end{align}
Summing up \eqref{ip1} and \eqref{ip2}, we obtain
\begin{align}
\int_S^T \int_0^1 (x\, \psi)_x  (F(\rho)+F(\xi)) \, dx\, dt =
\int_0^1 x\, \psi \, \left[F(\rho)\right]_S^T dx \notag \\ -\frac{1}{2}\int_S^T \int_0^1 a(x) x\psi (f(\xi)+f(\rho))
(\rho-\xi))\, dx\, dt.
\end{align}
Using the definition of $\psi$, we obtain
\begin{align}\label{eq-rho-xi1}
\int_S^T \int_{(0,1)\setminus Q_1}  (F(\rho)+F(\xi)) \, dx\, dt = -\int_S^T \int_{Q_1\cap (0,1)} (x\, \psi)_x
(F(\rho)+F(\xi)) \, dx\, dt +\int_0^1 x\, \psi \, \left[F(\xi)\right]_S^T dx \notag \\
 -\int_0^1 x\, \psi \, \left[F(\rho)\right]_S^T dx-\frac{1}{2}\int_S^T \int_0^1 a(x)x\psi (f(\xi)+f(\rho))(\rho-\xi)\, dx\, dt.
\end{align}
We now complete the expression of the energy $E_p$ in the left-hand side of the previous equality and,
since $\int_0^1 F(\xi)+F(\rho)\, dx=E_p$, it follows that
{\begin{align}\label{T123}
\int_S^T E_p(t) dt \leq& \underbrace{ \int_S^T \int_{Q_1\cap (0,1)} |\left(1-(x\, \psi)_x \right)| \left(F(\rho)+F(\xi)\right) \, dx\, dt}_{\mathbf{S_1}}+ \underbrace{\int_0^1 |x\,\psi | \,\left| \left[F(\rho)- F(\xi)\right]_S^T \right|dx}_{\mathbf{S_2}} \notag \\& \ \ \ \ \ \ \ \
\underbrace{+\frac{1}{2}\int_S^T \int_0^1 |a(x)x\psi|\left|(f(\rho)+f(\xi))\right||\rho-\xi|\, dx\, dt}_{\mathbf{S_3}} .
\end{align}}
We start by estimating $\mathbf{S_1}$. Since {$|\left(1-(x\,\psi)_x \right)| \leq C$}, we get
\begin{align}\label{T4}
\int_S^T \int_{Q_1\cap (0,1)} \left|1-(x\, \psi)_x \right| \left(F(\rho)+F(\xi)\right) \, dx\, dt
\leq C\int_S^T \int_{Q_1\cap (0,1)} \left(F(\rho)+F(\xi)\right) \, dx\, dt\leq C\mathbf{S_4},
\end{align}
where $\mathbf{S_4}$ has been defined in \eqref{T}.

As for $\mathbf{S_2}$, using the fact that $|x\psi|<1$ and the fact that $t\mapsto E_p(t)$ is non
increasing, one gets the following upper bound for $\mathbf{S_2}$,
\begin{align}\label{T1}
 \mathbf{S_2} \leq E_p(S)+E_p(T)\leq 2 E_p(S).
\end{align}
We finally estimate  $\mathbf{S_3}$. Recall that $q:=\frac{p}{p-1}$ denotes the conjugate exponent
of $p$. Using \eqref{young1} in Lemma~\ref{younglemma} with $A= a(x)\left|\rho-\xi\right| , B=\vert f(\xi)
\vert+ \vert f(\rho)\vert$ and $\eta=\eta_1$ where $\eta_1>0$ an arbitrary constant, it follows that
\begin{align}
\mathbf{S_3} &\leq C\int_S^T \int_0^1 a(x)\left|f(\rho)+f(\xi)\right| \left|\rho-\xi\right|\, dx\, dt \notag
\\ &\leq CC_p \eta_1^{q}\int_S^T \int_0^1 \left( F(\rho)+F(\xi)\right) \, dx\, dt + \frac{C C_p}{\eta_1^{p}} \int_S^T \int_0^1 a(x)|\rho -\xi|^p \, dx\, dt\notag
\\
&\leq CC_p\eta_1^q  \int_S^T E_p(t)dt +\frac{CC_p}{\eta_1^p} \int_S^T \int_0^1 a(x)|\rho -\xi|^p \, dx\, dt.
\label{termcomp}
\end{align}
Set $R= \max (|\rho|,|\xi|)$. Then, for every $0<\mu_1<1$, one has
\begin{align}\label{summu}
\int_S^T \int_0^1 a(x)|\rho -\xi|^p \, dx\, dt = \int_S^T \int_{|\rho -\xi|\geq R\mu_1} a(x)|\rho -\xi|^p \, dx\, dt + \int_S^T \int_{|\rho -\xi|< R\mu_1} a(x)|\rho -\xi|^p \, dx\, dt.
\end{align}
For the first integral term above, we have directly from Lemma \ref{lemmaA1} with
$a=\rho,\  b=\xi$ that
\begin{align}\label{geqmu}
\int_S^T \int_{|\rho-\xi |\geq R\mu_1} a(x)|\rho -\xi|^p \, dx\, dt \leq \frac{C_p}{\mu_1^{2-p}} \int_S^T(-E_p'(t))\, dt \leq \frac{C_p}{\mu_1^{2-p}} E(S).
\end{align}
As for the second integral term in \eqref{summu}, we have that
\begin{align}\label{leqmu}
\int_S^T \int_{|\rho-\xi |< R\mu_1} a(x)|\rho -\xi|^p \, dx\, dt &\leq \mu_1 ^p\int_S^T \int_{|\rho-\xi |< R\mu_1} a(x)R^p \, dx\, dt \notag \\
&\leq C_p\mu_1^p  \int_S^T \int_0^1( F(\rho)+F(\xi) ) dx dt \notag \\
& \leq C_p \mu_1^p \int_S^T E_p(t) dt.
\end{align}
Combining \eqref{summu}, \eqref{geqmu} and \eqref{leqmu}, we obtain that
\begin{align}\label{bornE_p'}
\int_S^T \int_0^1 a(x)|\rho -\xi|^p \, dx\, dt \leq C_p \mu_1^p  \int_S^T E_p(t) dt +\frac{C_p}{\mu_1^{2-p}} E_p(S).
\end{align}
By combining \eqref{bornE_p'} with \eqref{termcomp}, we obtain that
\begin{align}\label{T2}
\mathbf{S_3}\leq CC_p\left( \frac{\mu_1^p }{\eta_1^p} +\eta_1^q \right)\int_S^T E_p(t) dt +C\frac{C_p}{\eta_1^p \mu_1^{2-p}} E_p(S).
\end{align}
Gathering \eqref{T123}, \eqref{T4}, \eqref{T1} and \eqref{T2}, it follows that
\begin{align}
\int_S^T E(t)\,dt &\leq C\int_S^T \int_{Q_1\cap (0,1)} \left(F(\rho)+F(\xi)\right) \, dx\, dt\notag\\&+ C C_p\left( \frac{\mu_1^p }{\eta_1^p} +C\eta_1^q \right)\int_S^T E_p(t) dt + \left(C\frac{C_p}{\eta_1^p \mu_1^{2-p}}+2 \right)E_p(S).
\end{align}
We can choose $\eta_1>0$ and $\mu_1>0$ such that
\begin{equation}
CC_p\left( \frac{\mu_1^p }{\eta_1^p} +\eta_1^q \right)<\frac12,
\end{equation}
which proves \eqref{T}.
\begin{flushright}
\begin{small}
$\blacksquare$
\end{small}
\end{flushright}

\subsubsection{Second pair of multipliers}
The second set of multipliers given in Item $(m2)$ is used to handle the term $\mathbf{S_{4}}$ in \eqref{T} and it will lead us to the following lemma.
\begin{Lemma}
Under the hypotheses of Proposition~\ref{prop} with $\phi$ as defined in \eqref{psidef}, we have the
following estimate
{\begin{align}\label{T3}
\mathbf{S_{4}} \leq  C\frac{C_p}{\eta_2^p}\underbrace{\int_S^T\int_{Q_2 \cap
(0,1)} |z|^p \, dx\, dt}_{\mathbf{T_5}}+ C C_p\eta_2^q
\int_S^T E_p(t)\,dt +CC_pE_p(S),
\end{align}}
where $\eta_2$ is an arbitrary constant in $(0,1)$ and $C$ and $C_p$ are positive constants
whose dependence is specified in Remark~\ref{rem:Cp}.
\end{Lemma}
\emph{\textbf{Proof. }}
We multiply the first equation of \eqref{probw} by $\phi f'(\rho) z $, where $z$ is the solution of \eqref{prob} and we integrate over $[S,T]\times [0,1]$ to obtain
\begin{equation} \label{2ndmul}
\int_S^T\int_0^1 \phi f'(\rho) z (\rho_t - \rho_x
+\frac{1}{2}a(x)(\rho-\xi))\, dx\, dt=0.
\end{equation}
On one hand, we have that
\begin{align}\label{rhot}
\int_S^T\int_0^1\phi f'(\rho) z \rho_t \, dx\, dt =  \int_S^T\int_0^1 \phi \left(f(\rho)\right)_t z \, dx\, dt
=-\int_0^1 \int_S^T \phi f(\rho) z_t dtdx + \int_0^1 \phi \left[f(\rho)z \right]_S^T dx.
\end{align}
On the other hand, an integration by part with respect to $x$ yields
\begin{align}\label{rhox}
-\int_S^T\int_0^1  \phi f'(\rho) z \rho_x \, dx\, dt = - \int_S^T\int_0^1  \phi  z \left(f(\rho) \right)_x \, dx\, dt  =\int_S^T\int_0^1 (\phi z)_x f(\rho) \, dx\, dt,
\end{align}
and then
\begin{align}\label{that}
\int_S^T\int_0^1 (\phi z)_x\, f(\rho) \, dx\, dt = \int_S^T\int_0^1 \phi_x z f(\rho) \, dx\, dt + \int_S^T\int_0^1 \phi z_x f(\rho) \, dx\, dt.
\end{align}
Putting together \eqref{2ndmul},\eqref{rhot}, \eqref{rhox} and \eqref{that}
\begin{align}
- \int_S^T \int_0^1 \phi f(\rho) z_t \, dx\, dt + \int_0^1 \phi \left[f(\rho) z \right]_S^T dx +\int_S^T\int_0^1 \phi_x z f(\rho) \, dx\, dt \notag
\\
+ \int_S^T\int_0^1 \phi z_x f(\rho)\, dx\, dt +\frac{1}{2}\int_S^T\int_0^1 \phi f'(\rho) z a(x)(\rho - \xi) \, dx\, dt=0 .
\end{align}
Since
\begin{align}
\int_S^T&\int_0^1 \phi z_x f(\rho) \, dx\, dt - \int_S^T \int_0^1 \phi f(\rho) z_t \, dx\, dt=  \int_S^T\int_0^1 \phi \rho f(\rho) \, dx\, dt -2 \int_S^T\int_0^1 \phi z_t f(\rho) \, dx\, dt \notag \\
&= \int_S^T\int_0^1 \phi \rho f(\rho) \, dx\, dt -2 \int_S^T\int_0^1 \phi z_t f(\rho) \, dx\, dt,
\end{align}
it follows that
\begin{align}\label{rhomul2}
\int_S^T\int_0^1 \phi \rho f(\rho) \, dx\, dt= 2 \int_S^T\int_0^1 \phi z_t f(\rho) \, dx\, dt &- \int_0^1 \phi
\left[f(\rho) z \right]_S^T dx -\int_S^T\int_0^1 \phi_x z f(\rho) \, dx\, dt \notag \\
&-\frac{1}{2}\int_S^T\int_0^1 \phi f'(\rho) z\, a(x)(\rho - \xi) \, dx\, dt.
\end{align}
We proceed similarly after multiplying the second equation of \eqref{probw} by $\phi f'(\xi) z $ and,  following the same steps that led to \eqref{rhomul2}, we obtain that
\begin{align}\label{ximul2}
\int_S^T\int_0^1 \phi \xi f(\xi) \, dx\, dt=  \int_0^1 \phi \left[f(\xi) z\, dx  \right]_S^T &-\int_S^T\int_0^1
\phi_x z f(\xi) \, dx\, dt-2 \int_S^T\int_0^1 \phi z_t f(\xi) \, dx\, dt  \notag \\
&-\frac{1}{2}\int_S^T\int_0^1 \phi f'(\xi) z\, a(x)(\rho - \xi) \, dx\, dt .
\end{align}
We take the sum of \eqref{rhomul2} and \eqref{ximul2} and get
\begin{align}\label{mul2''}
\int_S^T\int_0^1 \phi \left(\rho f(\rho)+\xi f(\xi)\right) dx\, dt = -\int_S^T\int_0^1 \phi_x z \left(f(\rho) + f(\xi) \right)
\, dx\, dt-\left[ \int_0^1 \phi \left(f(\rho)-f(\xi)\right) z\, dx\right]_S^T  \notag \\
- \frac{1}{2}\int_S^T\int_0^1 \phi \left(f'(\rho) + f'(\xi) \right) z\, a(x)(\rho - \xi) \, dx\, dt  +2 \int_S^T\int_0^1 \phi z_t \left(f(\rho) - f(\xi)\right) \, dx\, dt .
\end{align}
Using the definition of  $\phi$ in \eqref{psidef} and the fact that $2z_t=\rho-\xi$, we derive that
{\begin{align}\label{mul2'}
\mathbf{S_{4}}  &\leq C\underbrace{\int_S^T\int_{Q_2 \cap (0,1)} \left|z \left(f(\rho) + f(\xi) \right)\right| \, dx\, dt}
_{\mathbf{T_1}}+ C_p\underbrace{\left|\left[ \int_0^1\left(f(\rho)-f(\xi)\right) z\, dx\right]_S^T\right| }
_{\mathbf{T_2}} \notag \\
&+C_p\underbrace{\int_S^T\int_{Q_2 \cap (0,1)} |\left(f'(\rho) + f'(\xi) \right) z\, a(x)(\rho - \xi)| \, dx\, dt}_{\mathbf{T_3}}  +C_p \underbrace{\int_S^T\int_0^1 \left|\phi (\rho-\xi) \left(f(\rho) - f(\xi)\right)\right| \, dx\, dt}_{\mathbf{T_4}} .
\end{align}}
We start by estimating ${\mathbf{T_1}}$. We have ${\mathbf{T_1}}\leq {\mathbf{T'_1}}$ where
\begin{equation}\label{eq:T'_1}
{\mathbf{T'_1}}:=  \int_S^T\int_{Q_2 \cap (0,1)} |z|(|f(\rho)|+|f(\xi)|) \, dx\, dt,
\end{equation}
which gives when using \eqref{young1} in Lemma~\ref{younglemma} with $A=|z|$ and
$B\in\{|f(\rho)|,|f(\xi)|\}$,
\begin{align}\label{S'}
&{\mathbf{T'_1}} \leq \frac{C_p}{\eta_2^p}\int_S^T\int_{Q_2 \cap (0,1)} F(z) \, dx\, dt + C_p
\eta_2^q\int_S^TE_p(t) dt,
\end{align}
where $\eta_2>0$ is arbitrary.
	
To estimate $\mathbf{T_2}$, we have by Young's inequality recalled in Lemma~\ref{lem:young} that
\begin{align}\label{S1est}
\left| \int_0^1 \phi \, \left(f(\rho)-f(\xi)\right) z dx \right|
&\leq  \int_0^1 |f(\rho)| |z| \, dx +  \int_0^1  |f(\xi)| |z| \, dx \notag\\
&\leq C_p\int_0^1\left( F(\rho)+F(\xi) \right) dx + C_p \int_0^1 |z|^p dx.
\end{align}
Using Poincar\'e's inequality,
\begin{align}\label{z^p}
\int_0^1 |z|^p dx\leq C\int_0^1 |z_x|^p dx &\leq C\int_0^1  |\rho+\xi|^p dx \leq CC_p\int_0^1\left(|\rho|^p + |\xi|^p\right) dx \leq CC_p \,E_p(t).
\end{align}
Combining \eqref{S1est} and \eqref{z^p} and the fact that $t\mapsto E_p(t)$ is non increasing, it follows
that
\begin{align}\label{S2}
\mathbf{T_2} \leq CC_pE_p(S).
\end{align}
As for $\mathbf{T_3}$, we first notice that for every $(\rho,\xi)\in\mathbb{R}^2$, one has
\begin{equation}
\vert \left(f'(\rho) + f'(\xi) \right) (\rho - \xi)\vert\leq C_p\left(\vert f(\rho)\vert+\vert f(\xi)\vert\right).
\end{equation}
It follows that $\mathbf{T_3}\leq C_p\mathbf{T'_1}$ which has been defined in \eqref{eq:T'_1} and which is upper bounded in \eqref{S'}.

Regarding $\mathbf{T_4}$, we use the fact that $\phi(x)\leq C a(x)$ for $x\in [0,1]$ to get
\begin{align}\label{S4}
\mathbf{T_4}\leq C\int_S^T\int_{\omega} a(x)(\rho - \xi) \left( f(\rho) -f(\xi)\right) \, dx\, dt
\leq C\int_S^T (-E'(t))dt\leq CE_p(S).
\end{align}
Combining \eqref{mul2'}, \eqref{S'}, \eqref{S2} and \eqref{S4},
the estimate \eqref{T3} is proved.
\begin{flushright}
\begin{small}
$\blacksquare$
\end{small}
\end{flushright}

\subsubsection{Third multiplier}
It remains to tackle the term $\mathbf{T_5}$ appearing in \eqref{T3}. To handle  it, we consider the
multiplier introduced in Item $(m3)$ and, in order to achieve future upper bounds, we will be needing
estimates of the $L^q$-norms of $v$ and $v_t$, where $q=\frac{p}{p-1}$, given in the following lemma.
\begin{Lemma}\label{elliptest}
For $v$ as defined in \eqref{ellipt}, we have the following estimates:
\begin{align}
&\int_0^1  |v|^{q}\, dx \leq CC_p E_p(t).  \label{estu}\\
\int_0^1  |v_t|^{q}\, dx &\leq C_p \left((p-2) \sigma E_p(t)+\frac{1}{\sigma^{p-2}}\int_0^1 \beta   |z_t|^p  dx\right) ,  \label{estut}
\end{align}
where $\sigma>0$ is an arbitrary positive constant  and $C$ and $C_p$ are
positive constants whose dependence is specified in Remark~\ref{rem:Cp}.
\end{Lemma}
\textbf{\emph{Proof. }} From the definition of $v$, one gets
\begin{equation}
v(t,x)=-x\int_x^1(1-s)\beta\, f(z)\, ds-(1-x)\int_0^xs\beta \, f(z)\, ds,\qquad x\in [0,1].
\end{equation}
One deduces that, by using H\"older's inequality,
\begin{equation}
\vert v(t,x)\vert^q\leq C\left(\int_0^1 \vert z\vert^{p-1}\, ds\right)^q\leq C\int_0^1 \vert z\vert^p\, ds,
\qquad x\in [0,1].
\end{equation}
Poincar\'e's inequality yields $\int_0^1 \vert z\vert^p\, ds\leq C\int_0^1 \vert z_x\vert^p\, ds$
and then \eqref{estu} after integrating over $x\in[0,1]$ and the definition of $E_p(t)$.

Similarly, one has
\begin{equation}
v_t(t,x)=-x\int_x^1(1-s)\beta\, z_tf'(z)\, ds-(1-x)\int_0^xs\beta \, z_tf'(z)\, ds,\qquad x\in [0,1].
\end{equation}
By using H\"older inequality and the fact that $\beta$ is bounded by $1$, one deduces that
\begin{equation}
\vert v_t(t,x)\vert^q\leq C_p\left(\int_0^1 \beta\vert z_t\vert\vert z\vert^{p-2}\, ds\right)^q
\leq C_p\int_0^1 \beta\vert z_t\vert^q\vert z\vert^{q(p-2)}\, ds,
\qquad x\in [0,1].
\end{equation}
If $p=2$, we have $q=2$ and get \eqref{estut} after integrating over $x\in[0,1]$. For $p>2$,
we apply Young's inequality with the pair of conjugate exponents $(p-1, \frac{p-1}{p-2})$ and conclude
as for \eqref{estu}.
\begin{flushright}
\begin{small}
$\blacksquare$
\end{small}
\end{flushright}
The next lemma shows the use of the third multiplier $v$.
\begin{Lemma}\label{lem:third}
Under the hypotheses of Proposition~\ref{prop}, with $v$ as defined in \eqref{ellipt}, we have the
following estimate,
\begin{align} \label{S1}
\underbrace{\int_S^T \int_{Q_2\cap(0,1)}|z|^p \, \, dx\, dt}_{\mathbf{T_5}} \leq  CC_p\left(\eta
\int_S^T E_p(t) \, dt+\frac1{\eta^r} E_p(S)\right),
\end{align}
where $r=\frac{2p^2-p-2}p$, $\eta$ is any real number in $(0,1)$ and $C$ and $C_p$ are positive constants
whose dependence is specified in Remark~\ref{rem:Cp}.
\end{Lemma}
\textbf{\emph{Proof.}}
We multiply the first equation of \eqref{probw} by $v$
\begin{equation} \label{3rdmul}
\int_S^T\int_0^1 v (\rho_t - \rho_x
+\frac{1}{2}a(x)(\rho-\xi))\, dx\, dt=0.
\end{equation}
First, an integration by part with respect to $t$ gives
\begin{align}\label{t1}
\int_0^1 \int_S^T v \rho_t\, dtdx= - \int_0^1 \int_S^T v_t \rho\, dtdx + \left[ \int_0^1 v \rho dx \right]_S^T.
\end{align}
Then, an integration by part with respect to $x$ yields
\begin{align}
-\int_0^1 v \rho_x \, dx = \int_0^1 v_x \rho \, dx =\int_0^1 v_x (z_x+z_t) \, dx =\int_0^1 v_x z_x \, dx+\int_0^1 v_x z_t\, dx.
\end{align}
We have that
\begin{align}
\int_0^1 v_x z_x \, dx= -\int_0^1 v_{xx} z \, dx = -\int_0^1 \beta |z|^p \, dx,
\end{align}
which gives that
\begin{align}\label{x1}
-\int_S^T \int_0^1 v \rho_x \, \, dx\, dt= -\int_S^T \int_0^1 \beta |z|^p \, \, dx\, dt+\int_S^T \int_0^1 v_x z_t\, \, dx\, dt.
\end{align}
Combining \eqref{3rdmul}, \eqref{t1} and \eqref{x1}, we obtain
\begin{align}\label{rho3}
\int_S^T \int_0^1 \beta |z|^p \, \, dx\, dt = \int_S^T \int_0^1 v_x z_t\, \, dx\, dt  - \int_0^1 \int_S^T v_t \rho\, dtdx + \left[ \int_0^1 v \rho dx \right]_S^T \notag
\\+ \frac{1}{2}\int_S^T \int_0^1 v a(x)(\rho-\xi)\, \, dx\, dt.
\end{align}
We next multiply the second equation of \eqref{probw} by $v$ and, following the same steps that
yielded \eqref{rho3}, we get
\begin{align}\label{xi3}
\int_S^T \int_0^1 \beta |z|^p \, \, dx\, dt = -\int_S^T \int_0^1 v_x z_t\, \, dx\, dt  + \int_0^1 \int_S^T v_t \xi\, dtdx - \left[ \int_0^1 v \xi dx \right]_S^T \notag
\\
+ \frac{1}{2}\int_S^T \int_0^1 v a(x)(\rho-\xi)\, \, dx\, dt.
\end{align}
Now taking the sum of \eqref{rho3} and \eqref{xi3}, we obtain
\begin{align}\label{sumrho3xi3}
2\int_S^T \int_0^1 \beta |z|^p \, \, dx\, dt = \int_0^1 \int_S^T v_t (\xi - \rho)\, dtdx - \left[ \int_0^1 v (\rho-\xi) dx \right]_S^T + \int_S^T \int_0^1 v a(x)(\rho-\xi)\, \, dx\, dt.
\end{align}
Using the definition of $\beta$, we obtain
\begin{align} \label{thirdmul}
2{\mathbf{T_5}} \leq  \underbrace{ \left|\left[ \int_0^1 v (\rho-\xi) dx \right]_S^T\right|}_{\mathbf{V_1}} + \underbrace{ \int_S^T \int_0^1 |v_t| |(\xi - \rho)|\, \, dx\, dt}_{\mathbf{V_2}} +\underbrace{ \int_S^T \int_0^1 |v a(x)(\rho-\xi)|\, \, dx\, dt}_{\mathbf{V_3}}.
\end{align}
We start by estimating $\mathbf{V_1}$. For fixed $t\in[S,T]$, we have, by using \eqref{estu}
\begin{equation}\label{eq:estV1}
\left| \int_0^1 v (\rho-\xi) dx \right| \leq \int_0^1 \left(|v| |\rho|+|v||\xi|\right) dx
\leq 2\frac{p-1}{p}\int_0^1 |v|^q dx +\int_0^1 \left(F(\rho)+F(\xi) \right)dx \leq CC_pE_p(t),
\end{equation}
and hence, since $E_p(T)\leq E_p(S)$, we get
\begin{align}\label{U1}
\mathbf{V_1} \leq \left| - \left[ \int_0^1 v(\rho-\xi) dx \right]_S^T \right| \leq CC_pE_p(S).
\end{align}
Using Young's inequality, we have for every $\eta>0$
\begin{align}\label{u2es}
\mathbf{V_2} \leq \int_S^T \int_0^1 \left( |v_t|
|\xi|+|v_t| |\rho| \right)\, dx\, dt
\leq  2\frac{p-1}{p \eta}\int_S^T \int_0^1 |v_t|^q \, dx\, dt + \eta\int_S^T \int_0^1  \left(F(\rho)+F(\xi)
\right)\, dx\, dt.
\end{align}
From \eqref{estut} and the fact that the definition of $\beta$ implies $\beta\leq C\, a$, we get for every $\sigma>0$ that
\begin{equation}
\int_S^T \int_0^1 |v_t|^q  \leq C_p\left((p-2)\sigma\int_S^TE_p(t)dt+\frac{C}{\sigma^{p-2}}\int_S^T \int_0^1a(x)|\rho-\xi|^pdx\, dt\right).
\end{equation}
Using \eqref{bornE_p'}, we obtain for every $\sigma,\mu_1>0$
\begin{equation}\label{otherhand}
\int_S^T \int_0^1 |v_t|^q  \leq CC_p\left[\left((p-2)\sigma+\frac{\mu_1^p }{\sigma^{p-2}}\right)\int_S^T
E_p(t) dt +\frac{1}{\mu_1\sigma^{p-2}} E_p(S)\right].
\end{equation}
Combining \eqref{u2es} and \eqref{otherhand}, we obtain for every $\eta,\sigma,\mu_1>0$
\begin{align}
\mathbf{V_2} \leq CC_p\left[\left((p-2)\frac{\sigma}{\eta}+\frac{\mu_1^p }{\sigma^{p-2}\eta}+\eta\right)
\int_S^T E_p(t) dt +\frac{1}{\mu_1\sigma^{p-2}\eta} E_p(S)\right].
\end{align}
Choosing $\sigma=\eta^2$ and $\mu_1=\eta^{2\frac{p-1}p}$, one gets, for every $\eta>0$
\begin{align}\label{U2}
\mathbf{V_2} \leq CC_p\left(\eta \int_S^T E_p(t) \, dt + \frac1{\eta^r} E_p(S)\right).
\end{align}
Finally, we estimate $\mathbf{V_3}$ in \eqref{thirdmul}. Using Young's inequality, we have, for every $\nu>0$,
\begin{align}
\mathbf{V_3}\leq CC_p\left(\nu \int_S^T \int_0^1 |v|^q \, \, dx\, dt + \frac1{\nu}\int_S^T \int_0^1  a(x)|\rho-\xi|^p\, \, dx\, dt\right),
\end{align}
which yields by using \eqref{estu} and \eqref{bornE_p'}, that for every $\nu,\mu>0$, one has
\begin{align}
\mathbf{V_3}
\leq CC_p\left((\nu + \frac{\mu^p}{\nu})\int_S^T E_p(t)dt  +\frac1{\nu\mu} E_p(S)\right).
\end{align}
Choosing $\mu^p=\nu$, one gets that for every $\eta>0$
\begin{align}\label{U3}
\mathbf{V_3}
\leq CC_p\left(\nu \int_S^T E_p(t)dt  +\frac1{\nu^{1+\frac2p}} E_p(S)\right).
\end{align}
Combining \eqref{thirdmul}, \eqref{U1}, \eqref{U2} and \eqref{U3} and taking $\nu=\eta<1$, we
obtain \eqref{S1}.
\begin{flushright}
\begin{small}
$\blacksquare$
\end{small}
\end{flushright}

\subsubsection{End of the proof of Proposition \ref{prop}}
\label{se:414}
Collecting \eqref{T}, \eqref{T3} and \eqref{S1}, we obtain for every positive $\eta_2,\eta_3$ and $\eta\in (0,1)$ that
\begin{equation}
\int_S^TE_p(t) dt \leq CC_p\left[\left(\eta_2^q+\frac{\eta}{\eta_2^p}\right)\int_S^T E_p(t) dt
+\left(1+\frac1{\eta_2^p\eta^r}\right)E_p(S)\right].
\end{equation}
Taking $\eta=\eta_2^{p+q}$ and fixing $\eta_2$ so that $2CC_p\eta_2^q=\frac12$,
we immediately get \eqref{expstab}. It is then standard to deduce that
there exists $\gamma_p > 0$ such that, for every $(z_0,z_1)\in X_p$, the energy $E_p$
associated with of the solution $z(t)$ of \eqref{prob} starting at $(z_0,z_1)$ satisfies the
following,
 \begin{equation}\label{expstabform}
E_p(t) \leq E_p(0) e^{1-\gamma_p t}, \qquad t\geq 0.
\end{equation}
That concludes the proof of Proposition \ref{prop}.

\subsection{Case where $\mathbf{1<p< 2}$}
The main issue to prove Theorem~\ref{thrmstab} in the case $p\in (1,2)$ (with respect to the case $p\in [2,\infty)$) is the trivial fact that
$p-2<0$ and hence the weights $f'(\rho)$ and $f'(\xi)$ used in the multipliers of Items $(m2)$ and
$(m3)$ may not be defined on sets of positive measure. As a consequence we cannot use these
multipliers directly and we have to modify the functions $f$ and $F$. This is why, we
consider, for $p\in (1,2)$, the functions $g$ and $G$ defined on $\mathbb{R}$, by
\begin{align}
g(y)&= (p-1)\int_0^y (|s| + 1)^{p-2} \, ds =\textrm{sgn}(y)\left[(|y|+1)^{p-1}-
1\right], \label{fepsdef}\\
G(y)&=\int_0^y g(s)\,ds =\frac{1}{p}\left[(|y|+1)^p-1\right]-|y|.
\label{Fepsdef}
\end{align}

It is clear that one has that
$\vert g(y)\vert\leq \vert f(y)\vert$ and $\vert G(y)\vert\leq \vert F(y)\vert$ for every $y\in\mathbb{R}$.
Finally, using the function $g$, we also modify
the energy $E_p$ by considering,  for every
$t\in \mathbb{R_+}$ and every solution of \eqref{prob}, the function ${\cal{E}}_p$ defined by
\begin{align}\label{Eepsdef}
{\cal{E}}_p(t)= \int_0^1 \left( G(\rho) +G(\xi) \right) \, dx.
\end{align}

We start with an extension of Proposition~\ref{prop1}.
\begin{Lemma}\label{convex}
For every $p\in (1,2)$,
\begin{itemize}
\item the function $g$ is an odd bijection from $\mathbb{R}$ to $\mathbb{R}$
with a continuous first derivative which is decreasing on $\mathbb{R}_+$;
\item the function $G$ is even, of class $C^2$ and strictly convex;
\item the energy $t\mapsto {\cal{E}}_p(t)$ is non-increasing on $\mathbb{R_+}$.
\end{itemize}
\end{Lemma}
\textbf{\emph{Proof.}} One has that, for every $x\in\mathbb{R}$,
\begin{align}
G''(x)=  g'(x)= (p-1)(|x|+1)^{p-2}>0.
\end{align}
It is clear that $g''$ is continuous and positive, which proves the strict convexity.
The last item follows after using Proposition~\ref{prop1} with ${\cal{F}}= \frac{G}{p}$ which admits a continuous
first derivative by what precedes.
\begin{small}
\begin{flushright}
$\blacksquare$
\end{flushright}
\end{small}
We define now the convex conjugate (cf.~\cite{Convexana2001}) of $G$ which we denote from
now on by $H$ and which is defined as the Legendre transform of $G$, i.e.,
\begin{align}\label{Gepsdef}
H(s)&:=\sup_{y\in\mathbb{R}}\{sy-G(y)\},\qquad s\in\mathbb{R}.
\end{align}
Since $G$ is of class $C^2$ with invertible first derivative $g$, one has that
\begin{align}\label{lgdt}
H(g(x))=\int_0^{g(x)}g^{-1}(s)ds=\int_0^x vg'(v)dv,
\end{align}
and
\begin{align}\label{lgdt1}
xg(x)= G(x) + H(g(x)), \quad x\in \mathbb{R}.
\end{align}
The second equality in \eqref{lgdt} is obtained using the change of variable $v=g^{-1}(s)$ and
\eqref{lgdt1} follows (for instance) by integration by part of the right-hand side of \eqref{lgdt1}.

The proof of Theorem~\ref{thrmstab} in the case $p\in (1,2)$ relies on the following
proposition which gives an estimate of the modified energy ${\cal{E}}_p$ of a strong solution and
which is similar to Proposition~\ref{prop}.
\begin{Proposition}\label{prop2}
Fix $ p \in (1,2)$ and suppose that Hypothesis $\mathbf{ (H_1)}$ is satisfied. Then there exist positive constants $C$ and $C_p$ such that, for every $(z_0,z_1) \in Y_p$ verifying
\begin{equation}\label{eq:assumption-epsilon}
E_p(0)\leq 1,
\end{equation}
we have the following energy estimate:
\begin{equation} \label{expstab2}
\forall \  0\leq S\leq T,\  \int_S^T {\cal{E}}_p(t)\, dt \leq C\,C_p{\cal{E}}_p(S).
\end{equation}
\end{Proposition}

We next develop an argument for Proposition~\ref{prop2}, which follows the lines of the proof of
Proposition~\ref{prop}. The main idea consists in replacing $f,F$ by $g,g$ and to control all the constants
$C_p$ involved in these estimates in terms of $p\in(1,\infty)$. We also provide a sketchy presentation where
we only precise details specific to the present case.

We fix $p\in (1,2)$ and $(z_0,z_1) \in Y_p$.
We recall that we have chosen $x_0=0$ as an observation point and let
$0<\epsilon_0 < \epsilon_1 < \epsilon_2 $ with the corresponding sets
$Q_i = ]1-\epsilon_i,1+\epsilon_i [$, $i=0,1,2$ as before.

As a consequence of \eqref{eq:assumption-epsilon} and Corollary~\ref{cor:decrease} and standard
estimates (such as the fact that ${\cal{E}}_p\leq E_p$), one deduces that
\begin{equation}\label{est-epsilon}
|z(t,x)|^p+{\cal{E}}_p(t)\leq C_p, \quad \forall t\geq 0,\ x\in [0,1],
\end{equation}
where $C_p$ is a positive constant that depends on $p$ only.
\subsubsection{{\textbf{First pair of multipliers} }}
For the first pair of multipliers, we change the function $f$ in Item $(m1)$ by the function
$g$ and hence  use $x\, \psi\, g(\rho),\ x\,\psi  g(\xi)$, where $\psi$ is defined
in \eqref{psidef}.
\begin{Lemma}
Under the hypotheses of Proposition~\ref{prop2}, we have the following estimate
\begin{align}\label{Tk2}
\int_S^T {\cal{E}}_p(t) dt \leq &\ C\,C_p {\cal{E}}_p (S) +
C \underbrace{\int_S^T \int_{Q_1\cap (0,1)}  \left(G
(\rho)+G(\xi)\right) \, dx\, dt}_{\mathbf{\overline{S}_4}}.
\end{align}
\end{Lemma}
\emph{\textbf{Proof.}}
Estimate~\ref{Tk2} is obtained by following the exact same steps as those given to derive \eqref{T}, with
the difference that we use the function $g$ instead of the function $f$.
By multiplying the first equation of \eqref{probw} by $x\, \psi\, g(\rho)$ and the second one
by $x\, \psi\, g(\xi)$, we perform the integrations by parts described to obtain \eqref{eq-rho-xi1}
with the function $f$ and, we are led to the similar equation
\begin{align}\label{ep1}
\int_S^T ( G(\rho)+ G(\xi))\, dt= \int_S^T \int_{Q_1\cap (0,1)} \left(1-(x\, \psi)_x \right)
(G(\rho)+ G(\xi)) \, dx\, dt +\int_0^1 x\, \psi \, \left[ G(\xi)- G(\rho)
\right]_S^T dx  \notag \\
 -\frac{1}{2}\int_S^T \int_0^1  a(x)x\psi( g(\xi)+g(\rho))(\rho-\xi)\, dx\, dt,
 \end{align}
 which yields that
 \begin{align}\label{T123c}
 \int_S^T {\cal{E}}_p(t) dt \leq &\underbrace{ \int_S^T \int_{Q_1\cap (0,1)} |\left(1-(x\, \psi)_x \right)|
 \left(G(\rho)+G(\xi)\right) \, dx\, dt}_{\mathbf{\overline{S}_1}}+
 \underbrace{\int_0^1 |x\,\psi | \,\left| \left[G(\xi)-G(\rho)\right]_S^T \right|dx}
 _{\mathbf{\overline{S}_2}}\notag\\
 &\underbrace{+\frac{1}{2}\int_S^T \int_0^1 |a(x)x\psi|\left| g(\xi)+ g(\rho)\right||\rho-\xi|\,
 dx\, dt}_{\mathbf{\overline{S}_3}} .
 \end{align}

 Using the fact that $\psi_x$ is bounded, we get at once that
 \begin{align}\label{St1}
\mathbf{\overline{S}_1} &\leq C\int_S^T \int_{Q_1\cap (0,1)}  \left(G
(\rho)+G(\xi)\right) \, dx\, dt\leq C\mathbf{\overline{S}_4},
\end{align}
where $\mathbf{\overline{S}_4}$ has been defined in \eqref{Tk2}.
Using now the fact that $ |x\psi|\leq 1$ and the fact that $t\mapsto {\cal{E}}_p(t)$ is non increasing, it
follows that
\begin{align}\label{St2}
\mathbf{\overline{S}_2} \leq {\cal{E}}_p(T)+{\cal{E}}_p (S) \leq 2{\cal{E}}_p (S).
\end{align}
As for $\mathbf{\overline{S}_3}$, we proceed as for the estimate of $\mathbf{S_3}$ by first
using \eqref{ineqafeps} and Lemma~\ref{lemmaA1gen} instead of Lemmas~\ref{younglemma} and ~\ref{lemmaA1} respectively.

In particular, we have the following estimate, which extends \eqref{bornE_p'} to the case $p\in(1,2)$ and
which holds for every $\mu_1\in (0,1)$,
\begin{equation}\label{eq:overbar-est}
\int_S^T \overbar{\cal{E}}_p(t)\, dt \leq C\, C_p \mu_1^p  \int_S^T {\cal{E}}_p(t) dt +C\, \frac{C_p}
{\mu_1^{2-p}} {\cal{E}}_p(S),
\end{equation}
where $\overbar{\cal{E}}_p(t)$ is defined by
\begin{equation}\label{eq:overbar}
\overbar{\cal{E}}_p(t)=\int_0^1a(x)g(z_t)\, dx.
\end{equation}

\subsubsection{{\textbf{Second pair of multipliers}}}
The goal of this subsection is to estimate $\mathbf{\overline{S}_4}$. To do so, we change the
function $f$ in Item $(m2)$ by the function $g$ and hence define the pair of multipliers:
$\phi g'(\rho)z,\phi g'(\xi) z$ where $\phi$ is defined in \eqref{psidef}.
\begin{Lemma}
Under the hypotheses of Proposition~\ref{prop2} and for $1<p<2$ with $\phi$ as defined in \eqref{psidef}, we have the following estimate:
\begin{align}\label{mul22}
\mathbf{\overline{S}_4}\leq	
C\,\frac{C_p}{\eta_2^p}\underbrace{\int_S^T\int_{Q_2 \cap
(0,1)} G(z) \, dx\, dt}_{\mathbf{\overline{T}_5}}+C\, C_p\eta_2^2
\int_S^T {\cal{E}}_p(t)\,dt +C\,C_p{\cal{E}}_p(S),
\end{align}
where $\eta_2$ is an arbitrary constant in $(0,1)$ and $C$ and $C_p$ are positive constants whose dependence are specified in Remark~\ref{rem:Cp}.
on $p$.
\end{Lemma}
\textbf{\emph{Proof.}}
Estimate~\eqref{mul22} is obtained by following the same steps as those given to derive \eqref{T3},
with $g$ instead of $f$.
By multiplying the first equation of \eqref{probw} by $\phi g'(\rho)z$ and the second one
by $\phi g'(\xi)z$, where $z$ is the solution of \eqref{prob}, we perform the integrations
by parts described to obtain \eqref{mul2''} with the function $f$ and we are led to the equation
\begin{align}\label{sum1}
\int_S^T\int_0^1 \phi \, \left(g(\rho)\rho+ g(\xi)\xi\right) dx\, dt &= - \ \int_S^T\int_0^1
\phi_x \,  z \, \left(g(\rho)+g(\xi)\right)\, dx\, dt +\left[ \int_0^1 \phi \, \left(g(\xi)-
g(\rho)\right) z dx \right]_S^T\notag \\
&\ \ \ \ -\frac{1}{2}\int_S^T\int_0^1 \phi (g'(\rho) +g'(\xi))z
a(x)(\rho - \xi) \, dx\, dt\notag \\
&+2\int_S^T\int_0^1 \phi \, \left(g(\rho)-g(\xi)\right) (\rho - \xi) dx\, dt.
\end{align}
According to \eqref{ineqxf}, one has
\begin{equation}
g(\rho)\rho+ g(\xi)\xi\geq g(\rho)+ g(\xi),\quad \forall (\rho,\xi)\in \mathbb{R}^2.
\end{equation}
Hence, also using the definition of $\phi$, it follows from \eqref{sum1} that
\begin{align}\label{sum2}
\mathbf{\overline{S}_4}&\leq
C\underbrace{ \ \int_S^T\int_{Q_2 \cap (0,1)} |z\left(g(\rho)+g(\xi)\right)|\, dx\, dt}
_{\mathbf{\overline{T}_1}} +\underbrace{\left|\left[ \int_{Q_2 \cap (0,1)}\left(g(\xi)-g(\rho)
\right) z dx \right]_S^T\right|}_{{\mathbf{\overline{T}_2}}}\notag\\
&+\frac{1}{2}\underbrace{\int_S^T\int_{Q_2 \cap (0,1)} \left(g'(\rho) + g'(\xi)\right)a(x)|
z(\rho - \xi)| \,  dx\, dt}_{\mathbf{\overline{T}_3}} +2\underbrace{\int_S^T\int_{Q_2 \cap (0,1)} |
\left(g(\rho)-g(\xi)\right) (\rho - \xi)| dx\, dt}_{{\mathbf{\overline{T}_4}}},
\end{align}
for some positive constant $C$. The above equation must be put in parallel with
\eqref{mul2'} where the term $\mathbf{\overline{T}_j}$, $1\leq j\leq 4$ in \eqref{sum2} corresponds to the
term $\mathbf{T_j}$ in \eqref{mul2'}.

The term $\mathbf{\overline{T}_1}$ is handled exactly as the term $\mathbf{T_1}$ while using
\eqref{ineqafeps} instead of Lemma~\ref{lem:young} in order to obtain
\begin{align}\label{eq:T1-epsilon}
&\mathbf{\overline{T}_1} \leq \frac{C_p}{\eta_2^p}\int_S^T\int_{Q_2 \cap (0,1)} g(z) \, dx\, dt +
C_p \eta_2^2\int_S^T{\cal{E}}_p(t) dt,
\end{align}
where $\eta_2>0$ is arbitrary.

We proceed similarly for the term $\mathbf{\overline{T}_2}$ by using \eqref{ineqafeps} instead
of Young's inequality and Corollary~\ref{cor:poincare-epsilon} instead of the standard Poincar\'e
inequality to obtain
{\begin{align}\label{eq:T2-epsilon}
&\mathbf{\overline{T}_2} \leq C_p{\cal{E}}_p(S).
\end{align}}

The term $\mathbf{\overline{T}_4}$ can also be treated identically as the term $\mathbf{T_4}$ to obtain
\begin{align}\label{eq:T4-epsilon}
&\mathbf{\overline{T}_4} \leq C{\cal{E}}_p(S).
\end{align}

We now turn to an estimate of $\mathbf{\overline{T}_3}$ which differs slightly from that of
$\mathbf{T_3}$ because of the appearance of the function $g'$. Using \eqref{est-epsilon} and
the second equation in \eqref{ineqx2}, one deduces that
\begin{equation}
\left(g'(\rho) + g'(\xi)\right)|z|\leq C_p|g(z)|,\quad t\in [S,T],\ x\in [0,1],
\end{equation}
where $C_p$ is a positive constant only depending on $p$. One derives that
\begin{equation}
\mathbf{\overline{T}_3}\leq C_p\int_S^T\int_{Q_2 \cap (0,1)} a(x)|g(z)|(|\rho|+|\xi|) \,  dx\, dt.
\end{equation}
Applying \eqref{ineqafeps} to the above, we end up with an estimate of
$\mathbf{\overline{T}_3}$ by exactly the right-hand side of \eqref{eq:T1-epsilon} and one concludes.
\begin{flushright}
\begin{small}
$\blacksquare$
\end{small}
\end{flushright}
{\subsubsection{{\textbf{Third multiplier} }}}
We finally turn to an estimation of the term $\mathbf{\overline{T}_5}$ and, relying on the
multiplier defined in Item $(m3)$, we get after changing the function $f$ by the function $g$
the multiplier (still denoted) $v$ solution of the following elliptic problem defined at every $t\geq 0$ by
\begin{equation}\label{ellipt'}
\left\lbrace
\begin{array}{ll}
v_{xx} = \beta g(z), & x\in [0,1] , \\
v(0)=v(1)=0, &
\end{array}
\right.
\end{equation}
where $\beta$ is defined in \eqref{psidef}.\\ \\
We will be needing the following estimates of $v$ and $v_t$ given in the next lemma.
\begin{Lemma}\label{elliptestc}
For $v$ as defined in \eqref{ellipt'}, we have the following estimates:
\begin{align}
&\int_0^1  H(v) dx \leq C C_p {\cal{E}}_p(t),  \label{estuc}\\
&\int_0^1  H(v_t)  dx \leq C C_p \overbar{\cal{E}}_p(t),  \label{estutc}
\end{align}
where $\overbar{\cal{E}}_p(t)$ is defined in \eqref{eq:overbar}.

\end{Lemma}
\textbf{\emph{Proof. }} From the definition of $v$, one gets
\begin{equation}
v(t,x)=-x\int_x^1(1-s)\beta\, g(z)\, ds-(1-x)\int_0^xs\beta \, g(z)\, ds,\qquad x\in [0,1].
\end{equation}
It immediately follows from the above that
\begin{equation}
\int_0^1  H(v) dx \leq H\left(C\int_0^1\beta|g(z)|\, ds\right)\leq
CC_pH\left(\int_0^1\beta|g(z)|\, ds\right),
\end{equation}
where we have used \eqref{ineqx2}whether $C_p\int_0^1\beta|g(z)|\, ds\geq M$ or not.
Since $H$ is (strictly) convex, one can apply Jensen's inequality to the right-hand side of the
above equation to get that
\begin{equation}
\int_0^1  H(v) dx \leq CC_p\int_0^1\beta H\left(g(z)\right)\, ds,
\end{equation}
and one derives \eqref{estuc} by using \eqref{ineqx2} together with \eqref{est-epsilon}.
Similarly, one has that
\begin{equation}
v_t(t,x)=-x\int_x^1(1-s)\beta\, z_tg'(z)\, ds-(1-x)\int_0^xs\beta \, z_tg'(z)\, ds,\qquad x\in [0,1],
\end{equation}
Upper bounding $|g'(z)|$ by $1$, one deduces that
\begin{equation}\label{eq:GGG}
\int_0^1  H(v_t)  dx \leq H(C\int_0^1\beta |z_t|)\leq
CC_pH(\int_0^1\beta |z_t|)\leq CC_pH\left(\int_0^1a(x) |z_t|\, dx\right),
\end{equation}
where we used the fact that $\beta(x)\leq C a(x)$ on $[0,1]$ and the convexity of $H$.

Since $\int_0^1a(x) |z_t|\, dx=\int_{|z_t|\leq M}+\int_{|z_t|> M}$,
we have according to \eqref{ineqx2}, \eqref{ineqxp} and H\"older's inequality that
\begin{align}
\int_0^1a(x) |z_t|\, dx&\leq C\,C_p\left(\int_0^1a(x)  g(|z_t|)\, dx+
\int_0^1a(x)g(|z_t|)^{\frac1{p}}\, dx\right)\notag\\
&\leq C\,C_p\,\left(\int_0^1a(x) g(|z_t|)\, dx+\overbar{\cal{E}}_p^{\frac1p}(t)\right).
\end{align}
By convexity of $H$, we obtain after plugging the previous equation into \eqref{eq:GGG} and
using Jensen's inequality that
\begin{equation}\label{}
\int_0^1  H(v_t)  dx \leq CC_p \left(\int_0^1a(x)H\left(g(|z_t|)\right)\, dx+
H(\overbar{\cal{E}}_p^{\frac1p}(t))\right).
\end{equation}
The first term in the right-hand side of the above inequality is clearly upper bounded by
$CC_p\overbar{\cal{E}}_p(t)$ thanks to \eqref{ineqG}. As for the second term, one has that
\begin{equation}\label{}
H(\overbar{\cal{E}}_p^{\frac1p}(t))\leq C C_pg(y_*),
\hbox{ where $y_*$ is defined by } g(y_*)=\overbar{\cal{E}}_p^{\frac1p}(t).
\end{equation}
By using \eqref{est-epsilon}, it follows that $g(y_*)\leq CC_p$ and elementary
computations using \eqref{ineqx2} and the fact that $\overbar{\cal{E}}_p\leq CC_p E_p$ yield that
$y_*\leq CC_p$. Hence $g(y_*)\geq CC_py_*$, i.e.,
$y_*\leq C C_p\overbar{\cal{E}}_p^{\frac1p}(t)$.  Since $g$ is convex and increasing on
$\mathbb{R}_+$, one gets after using  \eqref{ineqx2} that
\begin{equation}\label{}
H(\overbar{\cal{E}}_p^{\frac1p}(t))\leq C C_pg\left(C_p\overbar{\cal{E}}_p^{\frac1p}(t)\right)\leq C C_pg\left(\overbar{\cal{E}}_p^{\frac1p}(t)\right)\leq C C_p\overbar{\cal{E}}_p^{\frac2p}(t)\leq C C_p \overbar{\cal{E}}_p(t),
\end{equation}
where we have used repeatedly \eqref{est-epsilon} and \eqref{ineqx2}. This concludes the proof of \eqref{estutc}.
\begin{flushright}
\begin{small}
$\blacksquare$
\end{small}
\end{flushright}
We now use the multiplier $v$ in \eqref{probw} and we get the following result.
\begin{Lemma}
Under the hypotheses of Proposition~\ref{prop2} with $v$ as defined in \eqref{ellipt'}, we have the
following estimate
\begin{align} \label{S1'}
\mathbf{\overline{T}_5}\leq  C C_p\left(\eta^p
\int_S^T {\cal{E}}_p(t) \, dt+\frac1{\eta^s} {\cal{E}}_p(S)\right),
\end{align}
where $C_p$ is a positive number that depends on $p$ only
 and $\eta$ is any real number in $(0,1)$.
\end{Lemma}
\textbf{\emph{Proof.}} Proceeding as in the proof of Lemma~\ref{lem:third} to derive \eqref{thirdmul}, we
obtain
\begin{align} \label{thirdmul-epsilon}
2\int_S^T \int_{Q_2\cap(0,1)}zg(z)\ dx\ dt &\leq  \underbrace{ \left|\left[ \int_0^1 v (\rho-\xi) dx \right]_S^T\right|}_{\mathbf{\overline{V}_1}} + \underbrace{ \int_S^T \int_0^1 |v_t| |(\xi - \rho)|\, \, dx\, dt}_{\mathbf{\overline{V}_2}}\notag\\
&+\underbrace{ \int_S^T \int_0^1 |v a(x)(\rho-\xi)|\, \, dx\, dt}_{\mathbf{\overline{V}_3}}.
\end{align}
After using Fenchel's inequality \eqref{eq:fenchel} and \eqref{estuc}, one gets the following estimate for
$\mathbf{\overline{V}_1}$
\begin{equation}\label{eq:V1-ep}
\mathbf{\overline{V}_1}\leq CC_p {\cal{E}}_p(S).
\end{equation}
As for $\mathbf{\overline{V}_2}$, we first apply \eqref{ineqafeps1} (corresponding to the adaptation to the case $p\in (1,2)$ of the use of Young's inequality in \eqref{u2es}) to get that
$$
\mathbf{\overline{V}_2}\leq C\frac{C_p}{\eta^q}\int_S^T\int_0^1H(v_t)+CC_p\eta^p\int_S^T{\cal{E}}_p(t)dt,
$$
for every $0<\eta<1$. To handle the first integral term in the right-hand side of the above equation, we use \eqref{estutc} and \eqref{eq:overbar-est} to get that
$$
\mathbf{\overline{V}_2}\leq C\frac{C_p}{\eta^q\mu^{2-p}}{\cal{E}}_p(S)+CC_p\left(\frac{\mu^p}{\eta^q}+\eta^p\right)\int_S^T{\cal{E}}_p(t)dt,
$$
for every $0<\eta,\mu<1$. For $\mathbf{\overline{V}_3}$, we apply Fenchel's inequality, \eqref{estuc} and \eqref{eq:overbar-est} to get that
$$
\mathbf{\overline{V}_3}
\leq CC_p\left(\sigma^2+\frac{\lambda^p}{\sigma ^p}\right)\int_S^T{\cal{E}}_p(t)dt+C \frac{C_p}{\lambda^{2-p} \sigma^p}{\cal{E}}_p(S),
$$
for every $0<\lambda,\sigma<1$. One chooses appropriately $\lambda,\mu$ and $\sigma$ in terms of $\eta$ to easily conclude the proof of \eqref{S1'}.
\begin{flushright}
\begin{small}
$\blacksquare$
\end{small}
\end{flushright}
\subsubsection{End of the proofs of Proposition \ref{prop2} and Theorem~\ref{thrmstab} in the case $p\in (1,2)$}
It is immediate to derive \eqref{expstab2} by gathering \eqref{Tk2}, \eqref{mul22} and \eqref{S1'}
with a constant $C_p$ only depending on $p$. One deduces exponential decay of ${\cal{E}}_p$
exactly of the type \eqref{expstabform} with a constant $\gamma_p>0$ only depending on $p$ for
weak solutions verifying \eqref{eq:assumption-epsilon} for their initial conditions. Pick now any
$(z_0,z_1)\in X_p$ such that $E_p(0)=1$. One deduces that for every $t\geq 0$,
\begin{equation}\label{eq:fin-preuve1}
{\cal{E}}_p(t)\leq {\cal{E}}_p(0)e^{1-\gamma_p t}\leq e^{1-\gamma_p t},
\end{equation}
since ${\cal{E}}_p\leq E_p$. Set
$$
\lambda_p:=\left(\frac{p}8\right)^{\frac1p},
$$
and let $c_p$ be a positive constant such that
\begin{equation}\label{eq:int00}
G(x)>c_p F(x), \hbox{ if } \vert x\vert>\lambda_p.
\end{equation}
Note that such a constant $c_p>0$ exists according to the second equation in \eqref{ineqxp} and can be taken equal to $\frac{p-1}2$.

For every $t\geq 0$ and $x\in [0,1]$, let $R(t,x)=\max(\rho(t,x),\xi(t,x))$. It holds by elementary computations that
\begin{align}\label{eq:fin-preuve2}
\int_{R\leq \lambda_p}\left(F(\rho)+F(\xi)\right)dx&\leq \frac 14,\notag\\
\int_{R>\lambda_p}\left(F(\rho)+F(\xi)\right)dx&< \frac2{c_p}\int_{R> \lambda_p}\left(G(\rho)+G(\xi)\right)dx.
\end{align}
One deduces at once that 
\begin{equation}\label{eq:ouf00}
E_p(t)<\frac14+\frac2{c_p}{\cal{E}}_p(t),\quad \forall t\geq 0.
\end{equation}
Set 
$$
t_p:=\frac{\left(1+\ln(\frac8{c_p})\right)}{\gamma_p}.
$$
Then, using \eqref{eq:fin-preuve1}, it follows that $E_p(t)\leq \frac12$ if $t\geq t_p$ and then
$$
\Vert S_p(t_p)\Vert_{X_p}\leq \left(\frac12\right)^{\frac1p}<1,
$$
which implies that the $C^0$-semi-group $(S_p(t))_{t\geq 0}$ is exponentially stable for $p\in (1,2)$.
\begin{Remark}
From the argument, it is not difficult to see that $\gamma_p$ is bounded above and $c_p$ must tend to zero 
as $p$ tends to infinity. That yields that our estimate for $t_p$ tends to infinity as $p$ tends to one. Hence it is 
not obvious how to use our line of proof to get exponential stability for $p=1$.
\end{Remark}
\section{Case of a global constant damping}
Suppose now that we are dealing with a global constant damping, in other words $\omega=(0,1)$ and
\begin{align}
a(x)\equiv 2\alpha, \ \ \forall\ x \in (0,1),
\end{align}
where $\alpha$ is a positive constant. We then prove the following proposition.
\begin{Proposition}\label{prop:global} For $p=1$ or $p=\infty$, the semi-group $(S(t))_{t\geq 0}$ is
exponentially stable for a global constant damping if $\alpha\in (0,2)$.
\end{Proposition}
\textbf{\emph{Proof.}} For every $p\in (1,\infty)$, we perform a change of unknown function,
namely
$$
z(t,x)= e^{-\alpha t}v(t,x),\ \ x \in (0,1),\ \ t\geq 0,
$$
where $z$ is any solution of \eqref{prob} starting at $(z_0,z_1) \in X_p$. Clearly $v$ is a solution of
\begin{equation}\label{prob-v}
\left\lbrace
\begin{array}{lll}
v_{tt}-v_{xx}={\alpha ^2} v &\text{in }   \mathbb{R_+} \times (0,1), \\
v(t,0)= v(t,1)=0   &  t\geq 0 ,\\
v(0,\cdot)= z_0\ ,\ v_t(0,\cdot)=z_1+\alpha z_0.
\end{array}
\right.
\end{equation}
We use $E_p$ and $V_p$ to denote the $p$th-energies associated with $z$ and $v$ respectively. Since
$z_x=e^{-\alpha t}v_x$ and $z_t=e^{-\alpha t}(v_t-\alpha v)$, we get, after using Lemma~\ref{lem:Poincare}
and the following inequality (cf. \cite[Lemma 2.2]{adams2003})
\begin{equation}\label{eq:pp}
\vert a+b\vert^p\leq 2^{p-1}(\vert a\vert^p+\vert b\vert^p),\quad \forall (a,b)\in\mathbb{R}^2,
\end{equation}
that, for every $t\geq 0$,
\begin{equation}\label{eq:z-->v}
E_p(t)\leq e^{-\alpha pt}\left(2^{p-1}V_p(t)+2^p\alpha^p\int_0^1\vert v(x)\vert^p\, dx\right)\leq e^{-\alpha pt}\left(2^{p-1}+\frac{\alpha^p}{p^2}\right)V_p(t).
\end{equation}
On the other hand, for strong solutions of \eqref{prob-v}, one has after applying Corollary~\ref{cor:decrease}
to $v$ that, for every $t\geq 0$,
$$
V'_p(t)=\alpha^2\int_0^1v\left(\lfloor \rho \rceil^{p-1}-\lfloor \xi \rceil^{p-1} \right) dx,
$$
which yields, by using Young's inequality, that
\begin{align}
V'_p(t)\leq \alpha^2 \left( 2\frac{\eta^p}{p} \int_0^1  |v|^p dx+\frac{p}{q\eta^q} V_p(t)\right),
\end{align}
for every $\eta>0$. Using again Lemma~\ref{lem:Poincare} and the fact  that $v_x=\frac{1}{2}(\rho + \xi)$
together with \eqref{eq:pp}, we obtain that, for every $t\geq 0$,
$$
V'_p(t)\leq \alpha^2 \left(\eta^p K_p+\frac{p}{q\eta^q}\right) V_p(t),
$$
where we have set $K_p:=\frac1{p2^p}$.
The minimum with respect to $\eta$ of $ \eta^p K_p+\frac{p}{q\eta^q}$ is equal to $pK_p^\frac{1}{p}$, and one gets by using  Gronwall's lemma that
\begin{equation}\label{eq:Vp'}
V_p(t) \leq V_p(0) e^{\alpha^2 pK_p^\frac{1}{p} \, t}.
\end{equation}
Combining \eqref{eq:z-->v} and \eqref{eq:Vp'}, one gets that, for every $t\geq 0$,
$$
E_p(t)^{\frac1p}\leq (2+\alpha)^2e^{M_{\alpha}\, t}E_p(0)^{\frac1p},
$$
where
$$
M_\alpha:=-\alpha+\alpha^2 K_p^\frac{1}{p}=-\alpha\left(1-\frac{\alpha}{2p^{\frac1p}}\right).
$$
One concludes easily by letting $p$ tend either to one or $\infty$ and using an obvious density argument.
\begin{flushright}
\begin{small}
$\blacksquare$
\end{small}
\end{flushright}
 \bibliographystyle{plain}
\bibliography{biblio}
\appendix
\addcontentsline{toc}{section}{Appendix}
\section{Appendix}\label{appendix}
We next provide two classical inequalities.
\begin{Lemma}[\cite{Convexana2001}]\label{lem:young} \emph{(Young's inequality)}\label{younglemma}\\
Let $p>1$ and $q=\frac{p}{p-1}$ its conjugate exponent. Then, for every $A,B\in \mathbb{R}$ and
$\eta>0$, Young's inequality reads
\begin{align}
|AB|\leq \eta^p \frac{|A|^p}{p} + \frac{|B|^q}{ q\eta^q}. \label{young1}
\end{align}
\end{Lemma}
\begin{Lemma}[\cite{Convexana2001}]\label{fench} \emph{(Fenchel's inequality)}\label{lem:fenchel} \\
Let $a,b$ be two real numbers, and $f$ any function then it holds that
\begin{align}\label{eq:fenchel}
\left|a\,b\right| \leq  f(|a|)+ f^*(|b|),
\end{align}
where $f^*$ is the convex conjugate of $f$ defined by the Legendre transform as follows
\begin{align}
f^*(b)=sup_{x\in\mathbb{R}}\{b.x-f(x)\}, \qquad b\in \mathbb{R}.
\end{align}
Moreover, if $f$ is of class $C^2$, the derivative of $f^*$ is given by
\begin{align}
&\ \left[f^{*}\right]'(y)= \left[ f'\right]^{-1}(y),\qquad y\in\mathbb{R}.
\end{align}
\end{Lemma}
Note that Young' inequality is a particular instance of Fenchel's inequality, corresponding to the function $f(a)=\frac{\vert a\vert^p}p$ for $p>1$.

The next lemma states a technical result used several times in the paper.
\begin{Lemma} \label{lemmaA1} For $p>1$, there exists a positive constant $C_p$ such that, for every real numbers $a,b$
and $\mu\in (0,1)$ subject to $|a-b|\geq \max(|a|,|b|)\mu$, one has
\begin{align}\label{bornE'}
|a -b|^p\leq \frac{C_p}{\mu^{2-p}} (a-b)\left(f(a)-f(b) \right).
\end{align}	
\end{Lemma}
\textbf{\emph{Proof.}} With no loss of generality, we can assume that $\max(|a|,|b|)=|a|=R >|b|$ and have same sign.
Indeed, if $ab\leq 0$, then $|a -b|^p\leq 2^pR^p$ and $(a-b)(f(a)-f(b))\geq R^p$, hence
 \eqref{bornE'} is satisfied with $C_p\geq 2^{p})$.
Set then $h=1-\frac{ b}{a}$ and  $h\in (0,1)$.  Proving \eqref{bornE'} amounts to show that
there exists $C_p$ such that for every $h,\mu\in (0,1)$ with $h\geq  \mu$, it holds
\begin{align}\label{39}
h^{p-1} \leq \frac{C_p}{\mu^{2-p}} \left|1-(1-h)^{p-1} \right|.
\end{align}
Clearly the inequality holds for $h$ ``far away'' from zero for any $C_p$ large enough (w.r.t. one)
and hence it is enough to establish it for $h$ close to zero. By linearizing $(1-h)^{p-1}$, one must find
$C_p$ so that $h^{p-2} \leq \frac{C_p}{\mu^{2-p}}$ which indeed holds true.
\begin{flushright}
\begin{small}
$\blacksquare$
\end{small}
\end{flushright}
We can now state a lemma which is basic for our subsequent work.
\begin{Lemma}\label{lemmaFGxf}
Let $p\in (1,2)$. Then, the function $g,G$ and $H$ defined
in \eqref{fepsdef}, \eqref{Fepsdef} and \eqref{Gepsdef} satisfy the following relations:
\begin{description}
\item[$(i)$] for every $x\in\mathbb{R}$, one has
\begin{align}\label{eqG}
xg(x) = G(x) + H(g(x)).
\end{align}
\item[$(ii)$] for every $x\in\mathbb{R}$, it holds that
\begin{align}\label{ineqxf}
\frac{1}{2} \, x\,g(x) \leq  G(x) \leq \, x\, g(x).
\end{align}
\item[$(iii)$] There exists a positive constant $C_p$ only depending on $p$ such that,
for every $x\in \mathbb{R}$, one has
\begin{align}\label{ineqG}
C_p \, xg(x) \leq H(g(x)) \leq C_p xg(x).
\end{align}
\begin{align}\label{ineqGF}
C_p \, G(x) \leq H(g(x)) \leq C_p G(x).
\end{align}
\end{description}
\end{Lemma}
\noindent\textbf{Proof. }
Thanks to the parity properties of $g$ and $G$, it is enough to establish the several
relations only for $x\geq 0$.

Item $(i)$ is already proved in \eqref{lgdt1}. As for Item $(ii)$, the right inequality
\eqref{ineqxf} is immediate since $g$ is increasing. On the other hand, since $p<2$, we have
for all $0\leq s\leq x$ that
\begin{align}
g'(x)=(p-1)(x+1) ^{p-2} \leq (p-1)(s+1) ^{p-2}.
\end{align}
Integrating between $0$ and $x$, it follows that $xg'(x)\leq g(x)$, and then
$(xg)'(x)\leq 2g(x)$, which yields the left inequality \eqref{ineqxf} after an integration
between $0$ and $x$.
As for the proof of Item $(iii)$, it is clear that \eqref{ineqGF} follows from combining \eqref{ineqxf} and
\eqref{ineqG} and moreover, the right inequality in \eqref{ineqG} is an immediate consequence of
\eqref{eqG} since $g(x)\geq 0$ for $x\geq 0$. The proof for the left inequality in
\eqref{ineqG} is divided in two cases and can be deduced at once from the following estimates.
\begin{description}
\item[$(a)$] For every $M>0$ and real $x$ so that $|x|\leq M$, it holds
\begin{align}\label{ineqx2}
(p-1)(M+1)^{p-2}\frac{x^2}{2} \leq H(g(x)) \leq
(p-1)\frac{x^2}{2}. \notag \\
(p-1)(M+1)^{p-2}x^2 \leq x\, g(x)\leq (p-1)x^2 .\notag \\
(p-1)(M+1)^{p-2} \frac{x^2}{2} \leq  G(x)\leq (p-1)\frac{x^2}{2};
\end{align}
\item[$(b)$] for every $M>0$ such that $\left(1+\frac{1}{M}\right)^p<p$,
there exists a positive constant $C_p$ only depending on $p$ and $M$ so that, for every real $x$ verifying  $|x|>M$, one has
\begin{align}\label{ineqxp}
\left(\left(1+\frac{1}{M}\right)^{p-1}- \left(\frac{1}{M}\right)^{p-1} \right) |x|^p \leq x\, g(x)\leq |x|^p,\notag \\ \frac{1}{2}
\left(\left(1+\frac{1}{M}\right)^{p-1}- \left(\frac{1}{M}\right)^{p-1} \right) |x|^p \leq \, G(x)\leq  |x|^p,\notag \\
\left(1-\frac{1}{p}\left(1+\frac{1}{M}\right)^p\right)|x|^p\leq H(g(x)) \leq
\left(1+\frac{1}{M}\right)^{p-1}|x|^p.
\end{align}
\end{description}
\begin{Remark}
Note that the condition $\left(1+\frac{1}{M}\right)^p<p$ is only needed to get the third inequality of 
\eqref{ineqxp} only. Hence, the lower and upper bounds of $xg(x)$ and $G$ in \eqref{ineqxp} are valid for all 
$M>0$.
\end{Remark}
In turn, the set of inequalities in Item $(a)$ simply follows from the inequality
\begin{align*}
(M+1)^{p-2} \leq (s+1)^{p-2} \leq 1, \quad 0\leq s\leq x \leq M,
\end{align*}
and, after integrating between $0$ and $x$, by the use of the equations \eqref{fepsdef}, \eqref{Fepsdef}
and \eqref{lgdt}.

As for the set of inequalities in Item $(b)$, one first uses the explicit expressions of $xg(x)$ and
$G(x)$ given in \eqref{fepsdef}, \eqref{Fepsdef} to deduce that, for every $x\geq 0$,
\begin{align*}
x\,g(x)&=x^p \left[\left(1+\frac{1}{x}\right)^{p-1}-\left(\frac{1}{x}\right)^{p-1}\right],\\
G(x) &=x^p\left(\frac{1}{p}\left[\left(1+\frac{1}{x}\right)^p-\left(\frac{1}{x}\right)^p\right]-\left(\frac{1}{x}\right)^{p-1}\right).
\end{align*}
{Since $p<2$, the function $s\mapsto (1+s)^{p-1}-s^{p-1}$ is decreasing on $[0,\frac{1}M]$} and
then one gets the required bounds for $xg(x)$ in \eqref{ineqxp}. The upper and the lower bounds for
$g(x)$ in \eqref{ineqxp} are immediate and follow from combining the upper and the lower
bounds of \eqref{ineqxf} and the result above. Then the bounds for $H(g(x))$ in
\eqref{ineqxp} are simply obtained by combining the previous estimates with the relation
$H(g(x))=xg(x)-G(x)$.
\begin{flushright}
\begin{small}
$\blacksquare$
\end{small}
\end{flushright}
The next lemma is a particular instance of Fenchel's inequality which is used repeatedly in the paper.
\begin{Lemma}\label{lemmaafeps} For $p\in (1,2)$, there exist positive constants $C_p$ such that, for every $a,x\in \mathbb{R}$, it holds that
\begin{align}
|a\, x | &\leq  \frac{C_p}{\eta^p} G(a)+ \, C_p\,\eta^2 H(x),\label{ineqafeps0}\\
|a\, x | &\leq  C_p\eta^p G(a)+\frac{C_p}{\eta^q} H(x),\label{ineqafeps1}
\end{align}
and
\begin{align}\label{ineqafeps}
|a\, g(x) | \leq  \frac{C_p}{\eta^p} G(a)+ \, C_p\,\eta^2 g(x).
\end{align}
\end{Lemma}
\textbf{Proof :} Before proving the required inequalities, let us notice that one deduces from
\eqref{ineqx2} and \eqref{ineqxp} that there exists constants $C_p$ only depending on $p\in (1,2)$ such that
{\begin{align}
\hbox{if }\vert x\vert\leq C_p,& \hbox{ then }C_px^2\leq H(x)\leq C_p x^2,\label{eq:<Geps}\\
\hbox{if }\vert x\vert> C_p,& \hbox{ then }C_px^q\leq H(x)\leq C_p
x^q.\label{eq:>Geps}
\end{align}}
Since $g$ is convex, we
apply Fenchel's inequality given in \eqref{eq:fenchel} with $\frac{a}{\eta}$ and $\eta w$ with
$w\in\{x,g( x)\}$ to
obtain
\begin{align}\label{part3}
|a\, x| \leq G\left(\frac{a}{\eta}\right)+  H(\eta x).
\end{align}
Since both $H$ and $G$ are even functions, we assume with no loss of generality
that both $a$ and $x$ are non negative.

Using the estimates for $G$ and $H$ given in \eqref{ineqx2}, \eqref{ineqxp} and \eqref{eq:<Geps} and \eqref{eq:>Geps} respectively, we deduce that there
exists a positive constant $C_p$ only depending on $p\in (1,2)$ so that for every $a\geq 0$ and $\eta\in (0,1)$,
\begin{align}
G\left(\frac{a}{\eta}\right) &\leq C_p\max\left(\frac{1}{\eta^2},\frac{1}{\eta^p}\right) G(a)\leq \frac{C_p}{\eta^p} g(a),\notag\\
H(\eta x)&\leq C_p\max(\eta^2,\eta^q) H(x)\leq C_p\eta^2 H(x),
\label{eq:Gepsilon-eta}
\end{align}
and one immediately gets \eqref{ineqafeps0} from \eqref{part3} and \eqref{eq:Gepsilon-eta}. On the other hand, \eqref{ineqafeps} follows from \eqref{ineqafeps0} and \eqref{eq:Gepsilon-eta} after setting
setting $x=g(y)$ and using \eqref{ineqGF}.

Similarly, to get \eqref{ineqafeps1}, we start from
\begin{equation}\label{}
|a\, x| \leq G(\eta a)+  H\left(\frac{x}{\eta}\right),
\end{equation}
and we proceed as above to get the conclusion.
\begin{flushright}
\begin{small}
$\blacksquare$
\end{small}
\end{flushright}
As a corollary of the previous lemma, we have the following Poincar\'e-type of result.
\begin{Corollary}\label{cor:poincare-epsilon}
Let $p\in (1,2)$. Then there exists a positive constant $C_p$ such that, for every
absolutely continuous function $z:[0,1]\to\mathbb{R}$ so that $z(0)=0$, one has
\begin{equation}\label{eq:poincare-epsilon}
\int_0^1G(z(s))ds\leq C_p \int_0^1G(z'(s))ds.
\end{equation}
\end{Corollary}
\textbf{Proof.} With no loss of generality, we can assume that the right-hand side of \eqref{eq:poincare-epsilon} is finite. One has for every $x\in [0,1]$
\begin{equation}\label{}
G(z(x))=\int_0^x z'(s)g(z(s))ds.
\end{equation}
By applying \eqref{ineqafeps}, one gets that for every $x\in [0,1]$
\begin{equation}\label{}
G(z(x))\leq \frac{C_p}{\eta^p} \int_0^1 g(z'(s))ds+ C_p\,\eta \int_0^1
g(z(s))ds,
\end{equation}
for every $\eta>0$ and positive constants $C_p$ only depending on $p$. By integrating between $0$
and $1$ and then choosing appropriately $\eta$ one concludes.
\begin{flushright}
\begin{small}
$\blacksquare$
\end{small}
\end{flushright}

The following lemma is a useful extension of Lemma~\ref{lemmaA1} with $f,F$ replaced by
$g,g$.
\begin{Lemma} \label{lemmaA1gen}
For $p>1$, there exists a positive constant $C_p$ such that, for every real numbers $a,b$
and $\mu\in (0,1)$ subject to $|a-b|\geq \max(|a|,|b|)\mu$, one has
\begin{align}\label{ineqA1gen}
G(a-b)\leq \frac{C_p}{\mu^{2-p}}(a-b)\left(g(a)-g(b)\right)
\end{align}
\end{Lemma}
\textbf{Proof.}  Thanks to \eqref{ineqxf}, it i s enough to prove the
existence of $C_p>0$ so that
\begin{equation}\label{eq:lemA1gen1}
\vert g(a-b)\vert\leq \frac{C_p}{\mu^{2-p}}\vert g(a)-g(b)\vert,
\end{equation}
for every $a\geq b$, $\mu\in (0,1)$ such that $\vert a-b\vert \geq \mu R$ where
$R=\max(|a|,|b|)$. Assume first that $ab\leq 0$. Then the left-hand side of \eqref{eq:lemA1gen1} is
smaller than $g(2R)$ while $\vert g(a)-g(b)\vert\geq g(R)$. Clearly
$g(2R)\leq 2g(R)$ since $g$ is concave and hence \eqref{eq:lemA1gen1}
holds true in that case for any $C_p\geq 2$.

We next assume that  $a\geq b\geq 0$ and we consider $c=a-b$ instead of $b$. The assumption on
$a,b$ reads $c\geq \mu a$. Equation \eqref{eq:lemA1gen1} becomes
\begin{equation}\label{}
g(c)\leq  \frac{C_p}{\mu^{2-p}}\left(g(a)-g(a-c)\right).
\end{equation}
Note that the right-hand side of the above equation defines a decreasing function of $a$, once the other
parameters are fixed. It is therefore enough to consider the case $a=\frac{c}{\mu}$. By replacing $c$ by
$\frac{c}{\mu}$ in the explicit expression of $g$, we are led to prove the existence of $C_p>0$ so that
\begin{equation}\label{eq:lemA1gen2}
(\mu c+1)^{p-1}-1\leq \frac{C_p}{\mu^{2-p}} \left[(c+1)^{p-1}-\left((1-\mu) c+1\right)^{p-1}\right],
\end{equation}
for every $c>0$ and $\mu\in (0,1)$. By applying the mean value theorem to both sides of the above
equation and reordering the terms, \eqref{eq:lemA1gen2} reads
\begin{equation}\label{eq:lemA1gen3}
\left(\frac{\mu+\mu\eta_2}{\eta_1+1}\right)^{2-p}\leq C_p,
\end{equation}
for some $\eta_1\in (0,\mu c)$ and $\eta_2\in ((1-\mu) c,c)$ both depending on $c>0$ and $\mu$.
Assume first that $\mu c\leq 1$. Then clearly
\eqref{eq:lemA1gen2} holds true for any $C_p\geq 2^{2-p}$ according to \eqref{eq:lemA1gen3}. If now
$\mu c>1$, then $c>1$ and since the left-hand side of \eqref{eq:lemA1gen2} is smaller than
$(\mu c)^{p-1}$, we are left to find $C_p>0$ such that
\begin{equation}\label{}
\left(\frac{\mu+\mu\eta_2}{\mu c}\right)^{2-p}\leq C_p.
\end{equation}
The left-hand side of the above equation is again smaller than $2^{2-p}$ and one concludes.
\begin{flushright}
\begin{small}
$\blacksquare$
\end{small}
\end{flushright}
\begin{Lemma}\label{lem:Poincare}
Let $p\in(1,\infty)$. Then, for every $v \in W^{1,p}_0(0,1)$, it holds the
following Poincar\'{e} inequality
\begin{align}\label{Poincare}
\int_0^1 |v(x)|^p\,dx \leq  \frac1{p2^p}\int_0^1 |v'(x)|^p \, dx.
\end{align}
\end{Lemma}
\textbf{Proof.} For $x\in [0,\frac12]$, we have after using H\"{o}lder inequality,
\begin{equation}\label{}
|v(x)|^p=\left(\left| \int_0^x v'(s) \,ds\right|\right)^p \leq x^{\frac{p}{q}}\int_0^{\frac12} |v'(x)|^p\, dx.
\end{equation}
After integrating the previous between $0$ and $\frac12$, one gets
\begin{equation}\label{eq-012}
\int_0^{\frac12} |v(x)|^p\, dx\leq \frac1{p2^p}\int_0^{\frac12} |v'(x)|^p\, dx.
\end{equation}
For $x\in [\frac12,1]$, we start from $|v(x)|=\left| \int_1^x v'(s) \,ds\right|$ and get
\begin{equation}\label{eq-121}
\int_{\frac12}^1 |v(x)|^p\, dx\leq \frac1{p2^p}\int_{\frac12}^1 |v'(x)|^p\, dx,
\end{equation}
and finally \eqref{Poincare} by adding \eqref{eq-012} and \eqref{eq-121}.

\begin{flushright}
\begin{small}
$\blacksquare$
\end{small}
\end{flushright}

\end{document}